
\input amstex
\documentstyle{amsppt}

\NoRunningHeads
\magnification=1100
\baselineskip=29pt
\parskip 9pt
\pagewidth{5.2in}
\pageheight{7.2in}

\TagsOnLeft
\CenteredTagsOnSplits

\def\a{\Cal{A}}
\def\e{\Cal{E}}
\def\f{\Cal{F}}
\def\c{\Cal{C}}
\def\l{\Cal{L}}

\def\a{\Cal{A}}
\def\s{\Cal{S}}
\def\b{\Cal{B}}

\def\m{\Cal{M}}
\def\n{\Cal{N}}
\def\t{\Cal{T}}

\def\u{\Cal U}
\def\w{\Cal{W}}

\def\x{\Cal X}

\def\E{\Cal{E}}
\def\F{\Cal{F}}
\def\C{\Cal{C}}
\def\L{\Cal{L}}
\def\R{\Cal{R}}
\def\A{\Cal{A}}
\def\B{\Cal{B}}
\def\P{\Cal{P}}
\def\M{\Cal{M}}
\def\T{\Cal{T}}

\def\al{{\alpha}}
\let\alp=\al

\def\Bl{\bigl(}
\def\B{\Cal{B}}
\def\Bpd{\Cal{B}(P_d)\sta}
\let\BB=\B
\def\Br{\bigr)}

\def\br{\bigr)}
\let\bl=\Bl

\def\codim{\text{codim}}

\def\ctimej{{C\times J}}
\def\ctimew{{C\times W}}

\def\doubledual{^{\vee\vee}}

\def\etale{\`etale}
\def\endo{\Cal{E}nd^0}
\def\endpf{\hfill\qed\vskip8pt}
\def\etor{\e^t}

\def\Eb{\bar\e}
\let\eb=\Eb
\def\Ext{\text{Ext}}

\let\eps=\varepsilon

\def\dual{^{\vee}}
\def\d{\partial}

\def\ef{\e^f}
\def\fo{\f_1}
\def\ft{\f_2}

\def\germ{\text{germ}}

\def\hom{\Cal H\italic{om}}
\def\Hom{\text{Hom}}
\def\Hst{{$H$-stable}}
\def\Hs{{\Hst}}
\def\Hss{{$H$-semistable}}
\def\half{{1\over 2}}

\def\lra{\longrightarrow}
\def\ldel{_{\delta}}
\def\lalp{_{\alpha}}

\def\lep{_{\eps}}
\let\leps=\lep

\def\Lam{\Lambda}
\def\Loc{\Lambda_1^C}
\def\Ltc{\Lambda_2^C}

\def\len{\ell}
\def\lreg{_{\text{reg}}}

\def\JJ{\Cal I}

\def\mh{\!:\!}
\def\mapright#1{\,\smash{\mathop{\lra}\limits^{#1}}\,}
\def\mapto#1{\,\smash{\mathop{\to}\limits^{#1}}\,}

\def\mhid{\M_H(I,d)}
\def\mbid{\overline{\M}_H(I,d)}
\def\mhido{\mhid^0}

\def\M{{\frak M}}
\let\MM=\M

\def\MId{\M(I,d)}

\def\MId{\M(I,d)}
\def\mvid{\M(I,d)}
\def\mido{\M(I,d)^0}

\def\mdh{\M_d(H)}
\def\mb{\overline{\M}}
\def\mm{\Cal M_d}
\def\mmb{\overline{\Cal M}_d}
\def\mmm{\Cal M_{d-1}}
\def\mmmm{\Cal M_{d-2}}
\def\mmdo{\Cal M_{d-1}^0}
\def\mmo{\mm^0}
\def\mmmo{\mmm^0}

\let\Mido=\mido
\let\Mid=\mvid

\def\mbhid{\overline{\MM}_H(I,d)}

\def\mus{{$\mu$-stable}\ }
\def\muss{{$\mu$-semistable}\ }

\def\muxm{\mu_{XM}}

\def\N{\Cal{N}}

\def\nsr{\text{NS}_{\RR}}
\def\nsq{\text{NS}_{\QQ}}

\def\oplu#1#2{\smash{\mathop{\oplus}\limits^{#1}_{#2}}}
\def\OO{\Cal{O}}
\def\oppo{\OO_{\bold{P}}(1)}

\def\Pic{\text{Pic}}
\def\pixo{\pi^2_{X1}}
\def\pixt{\pi^2_{X2}}
\def\pito{\pi^2_1}
\def\pitm{\pi^2_{M}}
\def\pim{\pi^1_{M}}
\def\pix{\pi^1_{X}}
\def\pom{p_M}

\def\pfd{\Cal Z_1}
\def\ptfd{\Cal Z_2}
\def\Qu{\frak{Q}}

\def\pri{^{\prime}}
\def\PP{\bold P}

\def\rank{\text{rand }}

\def\resc{_{|C}}
\def\restc{_{|2C}}

\def\s{\Cal S}
\def\soo{\Cal S_1}
\def\so{\s_1}
\def\sox{\Cal S_1^x}

\def\st{\Cal S_2}
\def\sto{\Cal S_2^0}

\def\sta{^{\ast}}
\def\Sig{\Sigma}
\def\sig{\Sigma}

\def\supp{\operatorname{supp}}

\let\sub\subset
\let\ssub=\Subset
\def\Sl{\Cal{S}_l}

\def\Tor{\operatorname{Tor}}

\def\tilg{\tilde g}

\def\tilZ{{Z}}
\def\tenm{^{\otimes m}}

\def\tilZ{\tilde Z}

\def\tilE{{\tilde \E}}
\def\tile{{\tilde \e}}

\def\uptwo{^{\otimes 2}}
\def\upmo{^{-1}}
\def\uo{^0}
\def\util{\,\tilde{}\,}
\def\Vdel{V^{\del}}

\def\wedget{{\wedge^2}}
\def\xx{\chi}

\def\wtimej{{W\times J}}
\def\wtimec{{W\times C}}
\def\wb{\overline{\w}}
\def\wdel{\w^{\del}}
\def\wbdel{\overline{\w}^{\del}}

\def\zoreg{Z_{0,\text{reg}}}
\def\zreg{Z_{\text{reg}}}

\def\resc{_{|C}}
\let\inl=\iota
\def\inls{\inl_{\ast}}
\def\lam{\lambda}

\def\opmo{\OO_{\bold{P}}(-1)}

\def\lreg{_{\text{reg}}}
\def\del{\delta}

\def\dual{^{\vee}}

\def\wtj{_{W\times J}}

\def\ctj{_{C\times J}}
\def\tile{{\tilde \Cal{E}}}
\def\zpo{Z^{p_0}}

\def\ZZ{{\Bbb Z}}
\def\QQ{{\Bbb Q}}
\def\RR{{\Bbb R}}
\def\CC{{\Bbb C}}

\let\pro=\proclaim
\let\endpro=\endproclaim
\def\proof{{\it Proof.\ }}

\topmatter
\title
The first two Betti numbers of the moduli spaces of vector bundles
on surfaces
\endtitle

\author
Jun Li
\endauthor
\thanks This research was partially supported
by NSF grant DMS-9307892
\endthanks

\affil
Mathematics Department\\
University of California, Los Angeles
\endaffil

\email
jli@math.ucla.edu
\endemail

\endtopmatter
\document

\head 0. Introduction
\endhead

This paper is a continuation of our effort in understanding the geometry
of the moduli space of stable vector bundles. For any polarized
smooth projective surface $(X,H)$\ and for any choice of
$(I, d)\in \Pic(X)\times H^4(X,\ZZ)$, there is a coarse moduli space
$\M(I,d)\uo$\ of rank two $\mu$-stable (with respect to $H$)
locally free sheaves $\e$\ of $\wedget \e\cong I$\ and $c_2(\e)=d$.
This moduli space has been studied extensively recently.
One important discovery is that the moduli space
$\mvid^0$\ exhibits remarkable properties at stable range. To cite a
few, for arbitrary surface the moduli space $\M(I,d)^0$\ has the expected
dimension,
is smooth at general points and
is irreducible, and for a large class of surfaces of
general type $\mvid^0$\ are of general type,
all true for $d$\ sufficiently large [Fr, GL, Li2,
Do, Zh]. In this paper, we will investigate another aspect
of this moduli space. Namely, the Betti numbers
of $\Mido$. So far, there have been a lot of progress along this
direction based on two different approaches: Algebro-geometric
approach and gauge theoretic approach. The algebraic geometry
approach is relatively new. In [ES,Ki,Yo], they studied in detail the
Betti numbers of the moduli space of stable sheaves over
$\PP^2$\ (for the rank two and higher rank cases).
Beauville [Be] has a nice observation concerning some
rational surfaces and G\"ottsche and Huybrechts [GH] have worked out
the case for K3 surfaces.
The gauge theory approach has been around
for quite a while. To begin with, let $(M,g)$\ be a
compact oriented Riemannian four-manifold and let $P_d$\ be a
smooth SO(3) (or SU(2))-bundle
over $M$\ associated to a rank two vector bundle
of $c_1=I$\ and $c_2=d$. Consider the pair
$$\N(P_d)\sub\B(P_d)\sta, \tag 0.1
$$
where $\B(P_d)\sta$\ is the space of gauge equivalent classes of irreducible
connections on $P_d$\ and $\N(P_d)$\ is the
subspace of Anti-Self-Dual connections. By a celebrated
theorem of Donaldson, when $M=X$\ is an algebraic surface with a
Kahler metric associated to the ample divisor $H$,
$\Mido$\ is canonically isomorphic to $\N(P_d)$. The advantage of looking
at the pair (0.1) is that $H_{\ast}(\B(P_d)\sta)$\ is calculable,
at least modulo torsions,
in terms of the homotopy type of $X$\ and so does  $H_{\ast}(
\Mido)$\ if we know the induced homomorphism
$$
\ell(d)_i:H_i(\Mido,\ZZ)\lra H_i(\B(P_d)\sta,\ZZ).\tag 0.2
$$
In [AJ], Atiyah and Jones conjectured that for $M=S^4$\ and
SU(2)-bundle $P_d$, there is a sequence of (explicit) integers $\{q_k\}$\
such that for $d\geq q_k$, (0.2)
is an isomorphism for $i\leq k$. Later,  Taubes' work [Ta]
suggests that similar conjecture should hold for arbitrary
4-manifold with possibly different sequence
$\{q_k\}$. This conjecture has been confirmed for $S^4$, $\bold{CP}^2$\
and $K3$\ surfaces, see [HB, HBM$^2$, ES, GH, Ki, Ti1, Ti2, Yo].

In this paper, we will study $H_{\ast}(\Mido )$\ for arbitrary
algebraic surface. Due to technical difficulties, we are
unable to prove the generalized Atiyah-Jones conjecture for all
Betti numbers. Instead, we will prove the following theorems
that will determine the first two Betti numbers of the moduli space.

\pro{Theorem 0.1} For any smooth projective surface $(X,H)$\
and any $I\in \Pic(X)$, there is an integer $N$\ depending
on $(X,I,H)$\ so that whenever $d\geq N$, then
$$\ell(d)_i:H_i\bl\Mido;\QQ\br \lra H_i\bl \B(P_d)\sta;\QQ\br
$$
is an isomorphism for $i\leq 2$, where $P_d$\ is the SO(3) (or SU(2))-bundle
associated to the rank two complex vector bundle $E$\ with $\wedget E
\cong I$\ and $c_2(E)=d$.
\endpro

By calculating the first and second Betti numbers
of $\Bpd$, we get

\pro{Theorem 0.2}
With the notation as in theorem 0.1, then there is an integer $N$\
depending on $(X,I,H)$\ so that for $d\geq N$,
$\dim H_1\bl\Mido\br$\ and
$\dim H_2\bl\Mido\br$\ are $b_1$\ and $b_2+\half b_1(b_1-1)$\
respectively, where $b_i=\dim H_i(X)$.
\endpro

In algebraic geometry, there is a moduli space $\Mid$\ of
$H$-stable rank two sheaves $\e$\ with $\det \e\cong I$\
and $c_2(\e)=d$. $\Mid$\ is quasi-projective and
contains $\Mido$\ as a Zariski open subset.
We calculate the first two Betti numbers of $\Mid$\ as well.

\pro{Theorem 0.3}
With the notation as in theorem 0.1, then there is an integer $N$\
depending on $(X,I,H)$\ so that for $d\geq N$, $\dim H_1\bl\Mid\br$\
and $\dim H_2\bl\Mid\br$\ are $b_1$\ and $b_2+\half b_1(b_1-1)+1$\
respectively.
\endpro

There is a general principle [Mu] that explains why the Betti numbers
of $\Mid$\ take the form in theorem 0.3. For simplicity, let us assume
$\Mid$\ is projective and admits a universal family, say $\e$\ over $X\times
\Mid$. Then $\e$\ is expected to contain all information of $\Mid$. For
instance, the cohomology ring $H\sta(\Mid)$\ (with rational coefficient)
should be generated by the Kunneth components of $c_i(\e)$. Put
it differently, each $c_i(\e)$\ defines homomorphisms
$$\mu_i^{[\ast]}: H_{\ast}\bl X\br \lra H^{2i-\ast}\bl\Mid\br
$$
via slant product. Then the Mumford principle states that $H\sta(\Mid)$\ is
generated by the image of $\{\mu_i\}_{i\geq 2}$\ and within certain degree, say
$r(d)$\ ($r(d)\to\infty$\ when $d\to\infty$), their images obey no restraint
other than the obvious commutativity law of the cohomology ring $H\sta$.
In particular, if we look at $H^1(\Mid)$, then it should be generated (freely)
by
the images of $\mu_2^{[3]}\mh H_3(X)\to H^1(\Mid)$\ which
has dimension $b_1$\ by Poincare duality. For $H^2(\Mid)$, it
should be generated
freely by the wedge product of $H^1(\Mid)$, the image of $\mu_2^{[2]}\mh
H_2(X)\to H^2(\Mid)$\ and the image of $\mu_3^{[4]}\mh H_4(X)\to H^2(\Mid)$,
which span a linear space of total dimension $b_2+{1\over 2}b_1(b_1-1)+1$.

One motivation of the current work is to determine the Picard
group of the moduli space $\Mido$\ and $\Mid$. As is known,
$\Pic(\Mid)$\ is largely determined by $\dim H_1(\Mid)$\
and $\dim H_2(\Mid)$. In [Li3], we have determine
their Picard groups based on the information gained here.

Now we explain the strategy in establishing these Theorems.
In the following, we let $i=1$\ or $2$.
According to Taubes [Ta], for large $d$\ there are
canonical homomorphisms $\tau(d)_i$\ and $\tilde\tau(d)
_i$\ making the following diagram commutative:
\footnote{Taubes constructed a diagram using the space of based
connections. His diagram is identical to ours
because SO(3) and SU(2) are rational 3-sphere.}
$$\CD
H_i\bl\Mido\br @>{\ell(d)_i}>> H_i\bl\Bpd\br\\
@VV{\tau(d)_i}V @VV{\tilde\tau(d)_i}V\\
H_i\bl\M(I,d+1)^0\br @>\ell(d+1)_i>>
H_i\bl\BB(P_{d+1})\sta\br.\\
\endCD
\tag 0.3
$$
(Here and in the remainder of this paper, all homologies are with rational
coefficients unless otherwise is stated.)
Further, $\tilde\tau(d)_i$\ is an isomorphism, $\ell(d)_i$\ is
surjective for sufficiently large $d$\
and the composition of $\tau(\cdot)_i$'s
$$\tau(d,d+k)_i:
H_i(\Mido)\to H_i(\M(I,d+k)^0)
$$
has the property that
$$\ker\{\ell(d)_i\}\sub \ker\{\tau(d,d+k(d))_i\}
\tag 0.4
$$
for some $k(d)$. Thus, if
for large $d$\ the homomorphism $\tau(d)_i$\ is surjective,
$\tau(d)_i$\ must be an isomorphism for sufficiently large $d$.
Therefore by (0.4), $\ell(d)_i$\ will
be an isomorphism for sufficiently large $d$\ as well, thus
establishing theorem (0.1).

The homomorphism
$$\tau(d)_i: H_i\bl\Mido\br\lra H_i\bl\M(I,d+1)^0\br, \quad i\leq 2
\tag 0.5
$$
can be defined easily in our context.
We fix an $x\in X$\ and let $\sox\MM(I,d+1)\sub\MM(I,d+1)$\
be those $\e$\ such that $\e\doubledual/\e\cong\CC_x$.
$\sox\MM(I,d+1)$\ is a $\PP^1$-bundle over $\mvid$\
by sending $\e$\ to $\e\doubledual$. Let $V_0$\
be a general fiber.
Then the inclusion $V_0\sub\MM(I,d+1)$\ and the bundle
$\sox\MM(I,d+1)\to\MM(I,d)^0$\
induce a commutative diagram
$$\minCDarrowwidth{12pt}
\CD
0 @>>> H_i(V_0) @>>> H_i(\sox\MM(I,d+1)) @>>> H_i(\mido) @>>> 0 \\
@. @| @V r(d)_i VV @. @.\\
0 @>>> H_i(V_0) @>>> H_i(\MM(I,d+1)) @<<< H_i(\MM(I,d+1)^0).\\
\endCD
\tag 0.6
$$
When $d$\ is sufficiently large,
$H_i(\MM(I,d+1))$\ is a direct sum of the images of $H_i(V_0)$\
and $H_i(\MM(I,d+1)^0)$. Therefore, (0.6) induces a homomorphism
$H_i(\MM(I,d)^0)\to H_i(\MM(I,d+1)^0)$\ that is the mentioned
homomorphism $\tau(d)_i$. (See lemma 4.3 for details.)

{}From (0.6), $\tau(d)_i$\ is surjective if $r(d)_i$\ is surjective.
In this paper, we will prove the surjectivity of $r(d)_i$\ by
establishing the following two theorems:

\pro{Theorem 0.4}
With the notation as in theorem 0.1, then there is an integer $N$\
so that whenever $d\geq N$, then for $i=1,2$,
$$H_i\bl\mvid,\s\mvid\br=0,
$$
where $\s\mvid=\mvid-\mvid^0$.
\endpro

\pro{Theorem 0.5}
With the notation as in theorem 0.1, then there is an integer $N$\
so that whenever $d\geq N$, then for $i=1,2$\ and closed $x\in X$,
the image of
$$H_i\bl\s\mvid\br\lra H_i\bl\mvid\br
$$
is contained in the image of
$$H_i\bl\so^x\mvid\br\lra H_i\bl\mvid\br.
$$
\endpro

The strategy to establish theorem 0.4 is to apply the Lefschetz hyperplane
theorem to the moduli space.
The classical Lefschetz hyperplane theorem states that for any smooth,
projective
variety $Z$\ of complex dimension $n$\ and any smooth very ample divisor
$Z_1\sub Z$, the pair $(Z,Z_1)$\ has vanishing homology group up to degree
$n-1$. Concerning our situation, the ideal pair to look
at is $(\mvid, \s\mvid)$. But Lefschetz hyperplane theorem does
not apply directly to this pair because $\s\mvid$\ is definitely
not ample. Instead, we will first find an ample subvariety $\w$\ of
$\MId$\ and apply the generalized
Lefschetz hyperplane theorem to the pair $(\mido, \w)$\
to establish the surjectivity of
$$H_i\bl\w,\s\mvid\cap \w\br\lra
H_i\bl\mvid,\s\mvid\br,\quad i\leq 2.\tag 0.7
$$
The set $\w$\ has an explicit geometric description:
Let $C\in|nH|$\ be a fixed smooth divisor for some $n>0$.
Then $\w\sub\mvid$\ consists of those $\e$\ such that
$\e_{|C}$\ is not semistable.
By work of [Li1], there is a morphism
$$\varPsi_C:\mvid\lra \PP^M
\tag 0.8
$$
and a codimension $3g(C)-2$\ linear subspace $V\sub \PP^N$\
such that $\varPhi_C\upmo(V)=\w$. Hence we obtain the surjectivity of
(0.7) by applying the stratified Morse theory developed in [GM] to
$\varPsi_C$\ and $V\sub\PP^N$.

The next step is to show that
$$H_i\bl\w,\s\mvid\cap\w\br,\quad i\leq 2
\tag 0.9
$$
is trivial. The tactic is to construct explicitly a homology
between any class in (0.9) with the null class by
exploiting the fact that restriction to $C$\ of sheaves in $\w$\ are
not semistable. Here is an outline: For any locally free sheaf $\e\in\w$,
let $\l$\ be the destabilizing quotient sheaf of $\e_{|C}$\ and
let $\f$\ be the elementary transformation of $\e$\ defined by the
exact sequence
$$0\lra\f\lra\e\lra\l\lra0.
$$
Then $\e$\ can be reconstructed from $\f$\ via
$$0\lra\e\lra\f(C)\mapright{\al} I\otimes\L\upmo\lra 0.
$$
If we vary $\al$\ and $\l$, we get deformations of $\e$\ within $\w$.
In certain cases, we can deform $\e$\ to non-locally free sheaves this way.
This method was used by O'Grady in showing
that $H_0\bl\w,\s\mvid\cap\w\br=0$\ [OG1].
In this paper, we will work out this
construction in the relative case to prove the vanishing of (0.9).

A large portion of the current work is devoted to study
the singularities of various sets. This is necessary because
generalized Lefschetz hyperplane theorem only apply
to varieties with ``mild'' singularities.
In principle, the current approach should work for all homology
groups $H_i$\ through
a range that depends on the (local) topology of the singularities of $\w$.
For the moment, the author can only show that $\w$\ is locally irreducible
but has no idea whether it enjoys
any higher local connectivity.
Nevertheless, the local irreducibility of $\w$\ is
sufficient to show the vanishing of (0.9) and
thus establishing theorem 0.4. Theorem 0.5 is proved by carefully
studying the inclusion $\s\mvid\sub\mvid$.

The layout of the paper is as follows:
In \S 1, we will gather all relevant properties of the moduli space $\Mid$\
of which we will need. These include some discussion of singularities of
algebraic sets. In \S 2, by studying deformation of sheaves over curves, we
will show that the set $\w\sub\mvid$\
is locally irreducible. \S 3 is a refinement of [La, OG] in which we will
demonstrate how one can deform a family of
locally free sheaves to non-locally free sheaves and thus deriving the
vanishing of (0.9). The theorem 0.4 and 0.5
will be proved in \S 5. Most of
the materials concerning singularity are drawn from the excellent book of
Goresky and MacPherson [GM].

\head 1. Preliminaries
\endhead

In the first part of this section,
we will gather results of $\mvid$\
that are important to our study. Some of them have already appeared
or known to the experts and others are improvements of the
earlier results. We will give the reference to each result and
provide proof if necessary.
In the second part, we will review some materials concerning
topology of singular sets.

First, let us recall the convention that will be used throughout
this paper. In this paper, all schemes considered are of finite
type and are over complex number field. All points of schemes are
closed points. We will use Zariski topology throughout the paper
unless otherwise mentioned. Thus a closed subset is a union
of finite closed subvarieties. We will use algebraic subset to mean
finite union of locally closed subsets. By dimension of an algebraic set
we mean complex dimension.
We will only consider coherent sheaves in
this paper and will not
distinguish a vector bundle from the sheaf of its sections.

Throughout the rest of this paper, we fix a smooth algebraic surface
$X$\ and a line bundle
$I\in\Pic(X)$. Let $H$\ be an ample divisor on $X$. We say
a rank two sheaf $\E$\ is \Hst\ (resp. $H$-semistable)
if for any proper quotient sheaf $\e\to\f$, (i.e. it has non-trivial kernel)
we have
$$ {1\over \rank \e}
 \xx_{\E}(n) < {1\over \rank \f}\xx_{\f}(n) \quad \text{ (resp.
$\leq$)}
$$
for sufficiently large $n$, where $\xx_\E(n)=\xx(\E
\otimes H^{\otimes n})$\ is the value of the Hilbert polynomial of $\E$.
Note that $H$-semistable sheaves are necessarily torsion free.
Similarly, we say $\E$\ is $H$-\mus (resp. $H$-\muss)
if for any rank one {\sl torsion free} quotient sheaf  $\e\to \f$,
we have $\mu(\e)<\mu(\f)$\ (resp. $\leq $), where
$\mu(\E)={1\over \text{rank}\E}
c_1(\E)\cdot H$. We define stable and $\mu$-stable sheaves on curves
similarly. We also need the concept of $e$-stable.
For any constant
$e$, a rank two sheaf $\E$\ is said to be $e$-stable if for any rank one
{\sl torsion free} quotient sheaf $\e\to\f$, we have $\mu(\e)<\mu(\f)+e$.
One notices that for {\sl torsion free} sheaves,
$H$-$\mu$-stable implies \Hs\ and \Hss\ implies
$H$-$\mu$-semistable. This is not the
case for sheaves with torsion. For instance, it is easy to construct
a sheaf with torsion on curve that is $\mu$-stable but not stable.
In case the choice of $H$\ is apparent from the context,
we will simply call
them stable or $\mu$-stable. We agree that by unstable we mean
not semistable. According to [Gi], for any $d\in\ZZ$, there is a moduli scheme
$\mbid$\ of rank two $H$-semistable sheaves $\E$\ with $\det \E=I$\ and
$c_2(\E)=d$\ (modulo equivalence relation).
$\mbid$\ is projective. In the following, we will
freely refer a semistable sheaf $\e$\ as a point of $\mbid$.

There are several open subsets of $\mbid$\ that are relevant
to our study.
The first is the open subset $\mhid\sub\mbid$\ of all $H$-stable sheaves
(called the moduli of stable sheaves)
and $\mhid^0\sub\mhid$\ of all
$\mu$-stable locally free
sheaves. In most cases, $\mhido$\
is a Zariski dense open subset of $\mbid$. We let $\s\mbid=\mbid-\mhido$.
$\s\mbhid$\ contains non-locally free sheaves as well as locally
free but non $\mu$-stable sheaves.
For integer $l\geq 1$, we let $\s_l\mhid\sub \s\mhid$\ be the
set of non-locally free sheaves $\E$\ such that the length
$\len(\E^{\vee\vee}/\E)=l$, where $\E^{\vee\vee}$\ is the double dual
of $\E$, and let $\s^0_l\mhid\sub\s_l\mhid$\ be the subset of
those $\E$\ such that
$\E^{\vee\vee}/\E=\oplus_{i=1}^l \CC_{x_i}$\ with $x_i$\ distinct.
Note that $\s_l^0\mhid\sub \s_l\mhid$\ is open.

Usually, the algebraic subset $\s\mbhid\sub\mbhid$\
is not Cartier (Cartier means that set-theoretically
it is locally definable by one equation) which makes
the study of the topology difficult.
However, there is one situation that it does. Namely, when $H$\
is $(I,d)$-generic.

\pro{Definition 1.1}
1. An ample divisor $H$\ is called $(I,d)$-generic if for any
strictly $H$-semistable sheaf $\e$\ with $\det\e=I$\ and $c_2(\e)\leq
d+10$, $\e$\ is S-equivalent to a direct sum (of
rank one sheaves) $\l_1\oplus\l_2$\ such that $c_1(\l_1)=c_1(\l_2)\in
H^2(X,\RR)$.\hfil\break
2. Let $H_0$\ be any ample divisor. An ample divisor $H$\ is called
$(H_0,I,d)$-suitable if not only $H$\ is $(d,I)$-generic but
also has the property that any $H$-semistable sheaves $\e$\ with
$\det\e=I$\ and $c_2(\e)\leq d$\ are necessarily $H$-$\mu$-semistable.
\endpro

Here we face a dilemma: In proving the main theorems, we need to work
on $\mhid$\ inductively on $d$. Although for fixed
$(I,d)$\ there are plenty of $(I,d)$-generic ample divisors,
for fixed $H$, it will cease to be $(I,d)$-generic for large enough $d$,
assuming $\dim H^{1,1}(X)>1$. Thus we need to adjust
$H$\ constantly as $d$\ becoming large. To get by this,
we will work with a set of polarizations simultaneously.

To this end, some discussion on the selection of polarizations
is in order. First, because of
[Li1,p458], any two ample divisors $H_1$\
and $H_2$\ will give rise to (canonically) isomorphic moduli spaces
$\M_{H_1}(I,d)$\ and $\M_{H_2}(I,d)$\ if
$c_1(H_1)$\ and $c_1(H_2)$\ lie on the same (real) line in $H^2(X;\RR)$.
Now let
$$\nsr=\bigl( H^{1,1}(X,\RR)\cap H^2(X,\QQ)\bigr)\otimes_{\QQ}\RR,
$$
let $\nsr^+$\ be the ample cone
and let $\nsq$\ and $\nsq^+$\ be the intersection with $H^2(X,\QQ)$\
of the corresponding spaces. For any $\xi\in\nsq^+$, we define the
moduli space $\mb_{\xi}(I,d)$\ to be $\mbhid$\ for some ample
$H$\ such that $c_1(H)=n\xi$\ for some $n$.
By abuse of notation, in the following we will use $H\in\nsq^+$\
to mean $H$\ a $\QQ$-divisor with $c_1(H)\in\nsq^+$.
Next, let $H_0$\ be any ample divisor and let
$\c\leps\sub\nsq$\ be an $\eps$-ball in $\nsq$\ centered at $H_0\in
\nsq^+$,
after fixing an Euclidean metric on $\nsq$. For sufficiently small
$\eps>0$, the closure $cl(\c\leps)$\ of $\c\leps$\ in $\nsr$\ is
still contained in $\nsr^+$. We call such $\c\leps$\ precompact
neighborhood of $H_0\in \nsq^+$\ and denoted by $\c\leps\ssub\nsq^+$.

\pro{Lemma 1.2}
Let $H_0$\ be an ample line bundle and let $\c\ssub\nsq^+$\
be a precompact neighborhood of $H_0\in\nsq^+$.
Then for any choice of $(I,d)$, we can find an $(H_0,I,d)$-suitable
$\QQ$-ample divisor $H$\ in $\c$.
\endpro

\proof
It follows from theorem 1 on page 398 of [Qi] and the Hodge
index theorem.
\qed

{}From now on, we fix an $H_0\in\nsq^+$\ and a precompact neighborhood
$\c\ssub\nsq^+$\ of $H_0\in \nsq^+$. We will study moduli space
$\mhid$\ with $H$\ an arbitrary $\QQ$-divisor in $\c$\ and derive
estimate that depend on the set $\c\ssub\nsq^+$\ rather than individual
$H\in\c$.
We choose a smooth $C\in |2n_0H_0|$\ for some $n_0$\ such that
$$C\cdot C-|K_X\cdot C|- I\cdot C\geq 10.
\tag 1.1
$$
We will fix such a $C$\ and denote by $g$\ its genus in the remainder
of this paper.

Let $H\in\c$. Since usually the moduli space $\mhid$\ is singular, it
is convenient to work with a smooth subset of it:
$$\mm=\{\e\in\mhid\mid \e\ \text{is $\mu$-stable and}\
\Ext^2(\e\dual,\e\dual)^0=0\}.
\tag 1.2
$$
(For any sheaf $\E$, we let $\endo(\E)$\ be the sheaf of
traceless endomorphisms of $\Cal E$\ and let $\Ext^1(\E,\E)^0$\ be the
trace-less part of $\Ext^1(\E,\E)$.)
$\mm$\ is smooth. We let $\mm^0$, $\s\mm$\ and $\s_l\mm$\ etc. be
$\mm\cap\mhid^0$, $\mm\cap\s\mhid$\ and
$\mm\cap\s_l\mhid$\ etc. respectively.
For any constant $e<0$,
we let $\mbid_e\sub \mbid$\ be the set of all $e$-stable sheaves.
We summarize some properties of these sets in the following lemma.

\pro{Lemma 1.3}
There is an $N$\ depending on $(X,I,\c)$\ so
that whenever $d\geq N$, then for any $H\in\c$,
\roster
\item
$\mbhid$\ is normal, irreducible and has pure dimension $4d-I^2-
3\xx(\OO_X)$;
\item
$\mhid$\ is a local complete intersection scheme;
\item
Both $\mm$\ and $\so\mm$\ are smooth;
\item
$\mhido\sub \mbhid$, $\mhid\sub\mbhid$\ and $\mm\sub\mhid$\ are dense,
$\s_l^0\mhid$\ has dimension $\dim \mhid-l$\ and is dense in
$\s_l\mhid$\ for $l\leq 4$;
\item
The codimension
of the sets $\mbhid-\mhid_e$\ and $\mbhid-\mm$\
in $\mbhid$\ are at least $10g$, where $e=-2C^2$.
\endroster
\endpro

\proof
(1) and (2) were proved in [GL, Li2]. (3) is obvious and
(4) and (5) can be found in [Do, Fr, Li1, Qi, Zu].
\qed

We now introduce some subsets of $\mhid$\ associated $C$.
Let $2C\sub X$\ be the obvious non-reduced
subscheme supported on $C$. We define
$$ \Loc=\{\E\in\mhid\mid \Ext^2_X(\E\dual,\E^{\vee}(-2C))^0\ne0\}.
\tag 1.3
$$
$$\Ltc=\{\E\in\mhid\mid \E\resc\ \text{is locally free and}\
h^0(C,\endo(\E\resc))> g+1\}.
\tag 1.4
$$
$$\Lam_{\e}^C=\{ \f\in\mhid\mid \f\restc\cong \e_{|2C}\},\ \text{where $\e$\
is a sheaf locally free along $C$}.
\tag 1.5
$$

It is clear that all these sets are locally closed in $\mhid$.
Following [GL lemma 6.6],
$\Lam_{\e}^C$\ is a subscheme of $\mhid$.
Let $\f\in \Lam_{\e}^C$. The Zariski tangent space of
$\Lam_{\e}^C$\ at $\f$\ is isomorphic to the kernel of
the restriction homomorphism $\Ext^1_X(\f,\f)^0\to
\Ext_{2C}^1(\f\restc,\f\restc)^0$.

\pro{Lemma 1.4}
Let $\e$\ be any locally free sheaf on $2C$. Then
$\Lam_{\e}^C$\ is smooth and meets
$\so\mhid$\ transversally away from $\Loc$.
\endpro

\proof
This is a standard deformation calculation.
We shall leave the details to the readers.
\qed

Next we estimate the codimension of the sets $\Loc$\ and $\Ltc$.

\pro{Lemma 1.5}
There is an integer $N$\ depending on $(X,I,\c)$\ such that
whenever $d\geq N$\ and $H\in\c$, then
$\codim(\Loc,\mhid)\geq 10g$\ and
$\codim(\Ltc,\mhid)\geq 3g$.
\endpro

\proof
The proof given by [Do, Fr, Zu] can be adopted to cover the estimate $
\codim(\Loc,\mhid)\geq 10g$. Now we show that by choosing $N$\ large,
$$\codim(\Ltc,\mhid)\geq 3g.
$$
First of all, for any $\E\in\Ltc$, $\E\resc$\ is a direct sum of
line bundles, say $L_1\oplus L_2$\ with $\deg L_1>\deg L_2$\
and $h^0(L_1\otimes L_2^{-1})\geq g+3$\ by Riemann-Roch. Because $L_1\otimes
 L_2\cong I\resc$, the set $\{\E\resc\mid \E\in \Ltc\}$\ is
isomorphic to
an open subset of $\Pic(C)$. Next, for stable $\E\in\Ltc-\Loc$,
the tangent space of the set
$\{\E\pri \in\Ltc\mid \E\pri\resc\cong \E\resc\}$\
at $\E$\ is the kernel of $\Ext^1_X(\E,\E)^0\to\Ext^1_C(\E\resc,\E\resc)^0$\
which has dimension no more than
$$\dim\Ext^1_X(\E,\E)^0-\bl 3g-3+h^0(\endo(\E\resc))\br\leq
\dim\mhid-4g.
$$
Finally, since $\codim(\Loc,\mhid)\geq 3g$,
$$\codim(\Ltc,\mhid)\geq \min\{4g-\dim\Pic(C),3g\}=3g.
$$
This is exactly what we want.
\qed

In studying $\mhid$, we often need to use the local universal family. Let
$w\in\mhid$\ be any closed point. A local universal family of $\mhid$\ at $w$\
is an analytic (or \'etale)
neighborhood $g\mh U\to\mhid$\ of $w$\ and a flat family of sheaves
$\E_U$\ on $X\times U$\ flat over $U$\ such that for any $u\in U$,
the restriction sheaf $\E_{U|X\times \{u\}}$\ is represented
by the $g(u)\in\mhid$. By expressing $\mhid$\ as G.I.T. quotient of the
Grothendieck's quotient scheme and applying the \'etale slicing
theorem, we have

\pro{Lemma 1.6}
Any point $w\in\mhid$\ admits a local universal family.
\endpro

Now we discuss how to construct the morphism (0.8).
We fix an $N$\ given by lemma 1.3. For any $d\geq N$, let
$H\in\c$\ be $(H_0,I,d)$-suitable.
In [Li1], the author constructed a line bundle
$\l_C$\ over $\mbhid$\ and the associated morphism
$\varPsi_C:\mbhid\lra \PP^R$.
We summarize the property of this morphism as follows:

\pro{Lemma 1.7}\text{[Li1,\S2]}
Let $N$\ be given in lemma 1.3 and let $d\geq N$, $H\in\c$\
be $(H_0,I,d)$-suitable. Then
there is a line bundle $\l_C$\ on $\mbhid$\ of which the
following holds:
\roster
\item
For some large $m>0$, $H^0(\mbhid,\l_C\tenm)$\
is base point free and induces a morphism $\varPsi_C\mh\mbhid\to
\PP^R$;
\item
For any $\e\in\mbhid$, $\varPsi_C\upmo\bl\varPsi_C(\e)\br$\
consists of all $\f\in\mbhid$\ such that
$\f\dual\cong\e\dual$\ and $\ell\bl(
\f\doubledual/\f)\otimes\OO_x\br\cong \ell\bl(\e
\doubledual/\e)\otimes\OO_x\br$\ for all $ x\in X$;
\item
There is a codimension $3g-2$\ linear subspace $V\sub\PP^R$\
such that $\varPsi_C\upmo(V)$\ is exactly the set
$$\overline{\w}=\{\e\in\mbhid\mid \e_{|C}\ \text{is not semistable}\}.
$$
\endroster
\endpro

\proof
(1) and (2) follows directly from [Li1]. Now we prove (3).
Let $\M(C)$\ be the moduli space of semistable vector bundles over $C$
and let $\rho\mh \mbhid --\to\M(C)$\ be the rational map sending
$\e$\ to $\e_{|C}$\ when it is semistable. Then there is an ample
line bundle $L_C$\ on $\M(C)$\ such that $\rho\sta(L_C)\cong\l_C$\
over where $\rho$\ is defined. Further, for any $s\in H^0(\M(C),
L_C^{\otimes m})$, $\rho\sta(s)$\ extends to $\mbhid$\ with
vanishing locus
$$\{\e\in \mbhid\mid \e_{|C}\ \text{is not semistable or}\
\e_{|C}\in s^{-1}(0)\}.
$$
(See [Li1] for details.)
Thus, if we choose $3g-2$\ sections of $L_C^{\otimes m}$\
with no common vanishing locus, then the extensions of their pull
backs will define a codimension $3g-2$\ linear subspace $V\sub\PP^R$\
that has the desired property.
\qed

Note that in the proof, we only used the fact that all sheaves
in $\mbhid$\ are $H_0$-$\mu$-semistable. The extra requirement
that $H$\ be $(I,d)$-generic will be useful later
because of the following lemma.

\pro{Lemma 1.8}
Let $N$\ be as before and let $d\geq N$, $H\in\c$\ be
$(I,d)$-generic. Then the subset $\s\mbhid\sub\mbhid$\ is Cartier.
\endpro

\proof
The proof will appear in [Li3]. (See [DN] for a special case.)
\qed

In the remainder of this section, we will look more closely the local
geometry of various subsets of $\mbhid$.
This is necessary because later we are unable to
apply Lefschetz hyperplane theorem
directly to the space we want but rather a Zariski open
subset of it.

We still fix $H_0\in\c\ssub\nsq^+$, the $N$\ given by lemma 1.3 and
1.5. For any $d\geq N$, we choose an $(H_0,I,d)$-suitable $H\in\c$.
We let $\so\mhid\sub\mhid$\ and $\w\sub\mhid$\ be subsets defined
preceding definition 1.1 and in lemma 1.7 respectively.
Note that $\s\mbhid\sub\mbhid$\ is Cartier, $\mbhid$\ is normal,
$\mhid$\ locally is a complete intersection and $\w\sub\mhid$\
is defined by $3g-2$\ equations. We fix an embedding
$\mbhid\sub\PP^R$\ and an analytic Riemannian metric on $\PP^R$.
For any closed subset $A\sub\mbhid$\ and $\del>0$,
we define the $\del$-neighborhood of $A\sub\mbhid$, denoted by $A^{\del}$\
or $B\ldel(A)$, to be the set $\{z\in\mbhid\mid \text{dist}(z,A)<
\del\}$.

We first look at the subset $\so\mm\sub\mm$\ (cf. (1.2)).
$\so\mm$\ is a $\PP^1$-bundle over $X\times\mmdo$\
via the projection
$$\muxm:\e\in \so\mm \mapsto \bl\supp(\e\doubledual/\e),\e\doubledual\br.
$$
Let $V_0$\ be a general fiber of $\mu_{XM}$\ and let
$\n$\ be a normal slice of $\so\mm\sub\mm$\ along $V_0$.
(i.e. $\n$\ is an analytic neighborhood of $V_0\sub\tilde\n$, where
$\tilde\n\sub\mm$\ is a smooth analytic subvariety containing $V_0$\
that meets $\s\mm$\ transversally along $V_0$. For more details, see [GM].)

\pro{Lemma 1.9}
For $0<\del\ll 1$, $\n\cap B\ldel(V_0)$\ is a smooth analytic surface
that contains $V_0$\ as a (-2)-curve. In particular, $\n\cap B\ldel(V_0)-
V_0$\ is homeomorphic to $(\RR^4-\{0\})/\ZZ_2$.
\endpro

\proof
Since $V_0$\ is a general fiber of $\mu_{XM}$, $\mhid$\ is smooth along
$V_0$\ and for $0<\del\ll1$, $\n\cap B\ldel(V_0)$\ is a smooth
surface. The fact that $V_0\sub\n$\ is a (-2)-curve can be checked
directly by using the fact that sheaves in $V_0$\ are kernels of
$\e_0\to\CC_x$, where $V_0$\ lies over $(x,\e_0)\in X\times\M_H(I,d-1)^0$\
and the tangent bundle of $\mhid$\ at $z\in V_0$\ is
$\Ext^1(\e_z,\e_z)^0$. We shall omit the details here.
\qed

One type of  technical results that we need in the future says that
for closed subset $\Lam\sub\mbhid$\ of large codimension,
$$H_i\bl\mbhid-\Lam,\s\mbhid-\Lam\br\lra H_i\bl\mbhid,\s\mbhid\br
$$
is an isomorphism.  This type of results are certainly known to the
experts. Due to the lack of reference, we now provide proofs of them.

\pro{Lemma 1.10}
Let $Z\sub \PP^{n+r}$\ be any irreducible quasi-projective variety
of pure dimension $n$\ and $\Lam, V\sub Z$\ two closed algebraic subsets.
Assume $V$\ is Cartier, $\Lam$\ has codimension at least $k$\ and
$Z-V$\ is locally defined by at most $r+l$\ equations, then
$$H_i(Z-\Lam,V-\Lam)\lra H_i(Z,V)
$$
is an isomorphism for $i< (k-l)-1$\ and is surjective for $i=(k-l)-1$.
\endpro

\proof
Let $\s$\ be a Whitney stratification of $Z$\ by algebraic subsets so that
$V$\ and $\Lam$\ are union of strata. Let $S_1,\cdots,S_k$\ be strata of
$\Lam$\
with non-decreasing dimensions and let $\Lam_i=\cup_{j>i}S_j$. Then the
lemma follows if the homomorphism
$$H_i(Z-\Lam_j,V-\Lam_j)\lra H_i(Z-\Lam_{j-1},V-\Lam_{j-1})
$$
has the stated property for all $j$. Thus we only need to prove the
lemma under the assumption that $\Lam$\ is already a stratum in $\s$.
Since $\Lam\sub Z$\ is a stratum and $\Lam\sub Z$\ is closed,
there is a compact $\Lam_0\sub \Lam$\
such that $\bl(Z-\Lam)\cup \Lam_0,(V-\Lam)\cup\Lam_0\br$\
has the same homology group as $(Z,V)$. Hence
it suffices to show that
$$H_i(Z-\Lam,V-\Lam)\lra
H_i\bl(Z-\Lam)\cup \Lam_0,(V-\Lam)\cup\Lam_0\br
\tag 1.6
$$
has the desired property. Let $p_0\in \Lam_0$\ and let $N_{p_0}$\ be
the normal slice of $\Lam\sub Z$\ at $p_0$. We claim that the pair
$$\bl \partial B\leps(p_0)\cap N_{p_0},\partial B\leps(p_0)\cap N_{p_0}\cap
V\br
\tag 1.7
$$
is homologically $(k-l-2)$-connected for $0<\eps\ll 1$.
Indeed, the case when $p_0\in V$\ follows directly from theorem 2 and the
remark preceding it on page 156 of [GM]. When $p_0\not\in V$,
we let $(p_0,N_{p_0})\sub (0,\CC^R)$\
be an embedding and let $h_1,\cdots,h_t$,
$t=R-\dim N_{p_0}+r$, be the defining equations of $N_{p_0}$. Then the claim
follows from applying the same theorem to map $\pi\mh \CC^R\to\CC^{R+t}$,
$\pi(z)=(z,h_{\cdot}(z))$, and the linear subspace
$C^R\times\{0\}\sub\CC^{R+t}$.

Next by Thom's first isotopy lemma [GM], we can find a set
$U\sub (Z-\Lam)\cup\Lam_0$\ containing $\Lam_0$\ with a projection
$U\to\Lam_0$\ restricting to identity on $\Lam_0$\ such
that $U\to\Lam_0$\ is an $N_{p_0}$\ bundle over $\Lam_0$.
(i.e. $U$\ is a union of normal slices of $\Lambda\sub Z$\ at $p\in\Lam_0$.)
Because (1.7) is $(k-l-2)$-connected,
$(U-\Lam_0, (U-\Lam_0)\cap V)$\
is $(k-l-2)$-connected as well. Finally, we apply the Mayer-Vietoris sequence
to pairs $(Z-\Lam, V-\Lam)$\ and $(U, U\cap V)$. Because $(U,U\cap V)$\ has
trivial homology groups and $(U-\Lam_0, (U-\Lam_0)\cap V)$\ is
$(k-l-2)$-connected,
(1.6) is an isomorphism for $i< k-l-1$\ and surjective for $i=k-l-1$.
\qed

Let $\s$\ be a Whitney stratification of $Z$\ as before and let
$S\in\s$\ be a stratum. The following two lemmas concern the intersection
with $S$\ of representatives $\sig\sub Z$\ of elements in $H_i(Z)$.
We still denote by $N_p$\ a normal slice of $S\sub Z$\ at $p\in S$.

\pro{Lemma 1.11}
Let $S\in\s$\ be any stratum that is closed in $Z$. Then if
$\partial B\leps(p)\cap N_p$, $0<\eps\ll 1$, is $l$-connected (-1 connected
if it is non-empty and $-\infty$-connected if it is empty),
then any $v\in H_i(Z)$\ can be represented by a cycle $\sig\sub Z$\ such
that $\sig\cap S$\ has real dimension at most $i-l-2$.
\endpro

\pro{Lemma 1.12}
Let $S\in \s$\ be any stratum and let $q\in \overline{S}\sub Z$\ be
any point. Then if $\sig\sub Z$\ is a closed cycle contained in
$\cup_{i\in A}S_i$, where $A$\ is a subset of $\s$, and
$\sig\cap S$\ is discrete, then we can find a new representative $\sig'\sub
Z$\ of $[\sig]\in H_i(Z)$\ such that $\sig'\sub (\cup_{j\in A}S_j
-S)\cup\{q\}$.
\endpro

\proof
The proof of lemma 1.11 is obvious and we shall omit it. We now
prove lemma 1.12.
Let $U$\ be a cone neighborhood of $q\in Z$\
that respects the stratification $\s$. Let $p\in \sig\cap S$.
We choose a differentiable path $\rho\mh[0,1]\to S$\ connecting $p$\ with
a $p_1\in U$\ and let $N_{[0,1]}$\ be a continuous family of normal
slices of $S\sub
Z$\ along $\rho([0,1])$. By shrinking $N_{[0,1]}$\ if necessary, we can
assume $N_{[0,1]}$\ is homeomorphic to $N_p\times [0,1]$, say by
$\varPsi\mh N_p\times [0,1]\to N_{[0,1]}$. Next, by perturb $\sig$\ near
$p$, we can assume $\sig\cap B\leps(p)\sub N_p$\ for $0<\eps\ll1$. We fix a
sufficiently small $\eps>0$\ so that
$\varPsi((B\leps(p)\cap N_p)\times\{1\})\sub U$. Now let $A_1=\sig-B\leps(p)$,
$A_2=\varPsi(\Gamma\times[0,1])$, where $\Gamma=\partial(\sig\cap B\leps(p))$\
and
$A_3$\ is the cone over $\varPsi(\Gamma\times\{1\})$\ in $U$.
Then $\sig'=A_1\cup A_2\cup A_3$\ is a representative of $[\sig]$\ with
$\#(\sig'\cap S)<\#(S)$. By performing the above perturbation at each
$p\in \sig\cap S$, we will get a desired representative of
$[\sig]\in H_{\ast}(Z)$.
\qed

\head  2. Unstable sheaves over curves
\endhead

In this section, we will show that for any $\mu$-unstable
sheaves $\Eb$\ over $C$, the germ of the space
of $\mu$-unstable sheaves $\E$\ that are deformations
of $\Eb$\ is irreducible. To make it precise, we first recall
the notion of algebraic versal deformation space
of $\Eb$.

\pro{Definition 2.1}
Let $Z$\ be a projective scheme and let $\bar\E$\ be a rank two sheaf
over $Z$\ with $\det \bar\E\cong M$. An algebraic versal
deformation space of $\bar\E$\ (of fixed determinant) is a
tuple $(A,\F_A;0,\bar\E)$, where $0\in A$\ is
a quasi-projective scheme and $\F_A$\ is a flat algebraic family of
sheaves on $Z\times A$\ with
$\det \F_A\cong p_Z\sta M$, $p_Z\mh Z\times A\to Z$, such that
\roster
\item
the restriction of $\F_A$\ to $Z\times\{0\}$, say $\F_0$, is isomorphic
to $\bar\E$, and further the Kodaira-Spencer map $T_0A\to
\Ext^1(\F_0,\F_0)^0$\ induced by the family $\F_A$\ is an isomorphism;
\item
For any marked variety $s_0\in S$\ coupled with a
flat family of sheaves $\E_S$\ on $Z\times S$\ with
$\det \E_S\cong p_Z\sta M$\
and $\E_{s_0}\cong \bar\E$,
there is an analytic neighborhood $U$\ of $s_0\in S$\ and an analytic
map $\eta\mh (U,s_0)\to (A,0)$\ so that the restriction of $\E_S$\ to
$Z\times U$\ is isomorphic to $(1_Z\times \eta)\sta \F_A$,
extending the given isomorphism $\F_0\cong \bar\E$\ and
$\E_{s_0}\cong \bar\E$.
\endroster
\endpro

\pro{Lemma  2.2}
Assume $Z$\ is a projective curve and $H^0(Z,\OO_Z)=\CC$. Then for
any rank 2 sheaf $\E$\ on $Z$\ of $\det\E\cong M$,
there is an algebraic versal deformation space
(of fixed determinant) of $\E$.\endpro

\proof
The existence of this space in our setting is known for long
time (see [At]). Here we outline the proof since we need some
properties of this deformation space that can not be found in
reference.

By choosing a sufficiently ample line bundle $H$\ on $Z$, we
can express $\eb$\ as a quotient sheaf of  $\R=\oplus^N H^{-1}$\
with $N=h^0(\eb\otimes H)$.
Let $\Qu$\ be the Grothendieck's Quot-scheme parameterizing
all quotient
sheaves $\F$\ of $\R$\ with $\det \F\cong M$.
We fix a point $0\in\Qu$\ so that $\f_0\cong \eb$\ and
the associated quotient sheaf
$\R\mapright{\sigma_0}\f_0$\ induces isomorphism $\CC^N\cong
H^0(\f_0\otimes H)$.
Because $\dim Z=1$\ and $H$\ is sufficiently ample,
$\Qu$\ is smooth at $0$\ [Ma, p594]
and the tangent space of $\Qu$\ at $0$\ belongs to the exact sequence
$$0\to \Hom(\f_0,\f_0)\to \Hom(\R,\f_0)\mapright{j} T_0\Qu\to
\Ext^1(\f_0,\f_0)^0\to 0.
\tag 2.1
$$
Now we let $A\sub\Qu$\ be an affine smooth subvariety
containing $0$\
so that the induced homomorphism
$$T_0A\hookrightarrow T_0\Qu\lra \Ext^1(\f_0,\f_0)^0
$$
is  an isomorphism. Let $\F_A$\ be the restriction to
$Z\times A$\ of the universal quotient family.
Then the data $(A,\F_A;0,\F_0)$\ satisfies 1) of the definition 2.1.

Now, we show that they also satisfy 2) of the definition 2.1.
Let $s_0\in S$\ and $\E_S$\ be given in (2).
Because $H$\ is sufficiently ample,
by shrinking $S$\ is necessary, there is a morphism
$g\mh S\to\Qu$\ with $g(s_0)=0$\ so that $\E_S$\ is isomorphic to
the pullback of the universal quotient family.
It remains to find a neighborhood $U$\ of $0\in\Qu$\ and
an analytic map $\pi_{+}\mh U\to A$\
so that $\pi_+\circ g$\ is the desired map $\eta$. We
argue as follows: $\Qu$\ is a $G$-scheme with
$G=GL(N,\CC)=\Hom(\R,\R)$\ (because $H^0(Z,\OO_Z)=\CC$).
Let $G_0$\ be the stabilizer of $0\in\Qu$\
and let $N$\ be a normal slice of $G_0\sub G$\ at $1\in G_0$.
Then the orbit $N\cdot 0$\ is smooth at $0$\ and its
tangent space at $0$\ is the image $j(\Hom(\R,\f_0))\sub T_0\Qu$.
Thus $N\cdot 0$\ and $A$\ meets transversally at $0$.
In particular, there is an analytic neighborhood $V_-$\ of $1\in N$\
and an analytic neighborhood $V_+$\ of $0\in A$\ so that $V_-\times V_+
\to\Qu$, $(u_-,u_+) \mapsto u_-\cdot u_+$,
is one-to-one. Let $\u=\text{image of}\ V_-\times V_+$\ and let
$$\pi_{\pm}: \u \lra V_{\pm} \tag 2.2
$$
be the induced projections.
Then for any analytic variety $B$\
and any analytic map $\xi\mh  B\to V_{-}\times V_{+}\sub\Qu$,
$\xi(z)=f(z)\cdot(\pi_+\circ\xi(z))$,
where $f(z)=\pi_-\circ\xi(z)\in N$. Hence the sheaf $\e_B$\ in
the pullback quotient family
$h_B\mh p_Z\sta \R\to \E_B$\ via $1_Z\times \xi$\
and the sheaf $\e_B\pri$\ in the pullback quotient family
$h_B\pri\mh p_Z\sta\R\to\E_B\pri$\ via
$1_Z\times( \pi_+\circ \xi)$\ are isomorphic (analytically).
Back to $g\mh S\to\Qu$, we let $U=g\upmo(\u)$\ and
$\eta=\pi_+\circ g\mh U\to A$. The previous reasoning shows that
$(\eta,U)$\ is the desired map.
\qed

\pro{Remark 2.3}
For $w\in V_+\sub A$, because
$T_w\Qu\to\Ext^1(\F_w,\F_w)^0$\ is surjective and
$T_w( V_-\cdot w)\to \Ext^1(\F_w,\F_w)^0$\ is trivial,
$T_wA\to\Ext^1(\F_w,\F_w)^0$\ must be surjective.
\endpro

{}From now on, we fix a smooth connected curve $C$\ of genus $\geq4$,
a rank two $\mu$-unstable sheaf $\bar\E$\ over $C$\
of even degree
and a versal deformation space $(A,\F_A;0,\bar\E)$\ of $\bar\E$.
By tensoring $\bar\e$\ with $M\upmo$, where $M^{\otimes 2}\cong
\det\bar\e$, we can assume without loss of generality that
$\det\bar\e\cong\OO_C$.
We first stratify $A$\ as follows: For any $l>0$,
we let $\s_l$\ be the set of all rank two $\mu$-unstable sheaves
with $\det \E\cong\OO_C$\ whose destabilizing quotient sheaves
have degree $-l$\ and let
$\s_l^0$\ be those in $\s_l$\ that have $H^0(\endo(\E))=0$.
In light of our discussion in \S 1, we are interested in the set
of points in $A$\ whose associated sheaves are $\mu$-unstable. Namely,
$$A_0=\{w\in A\mid \F_w\ \hbox{is $\mu$-unstable}\}.\tag 2.3
$$
$A_0$\ is a locally closed algebraic set. Since $\bar\e$\ is $\mu$-unstable,
$0\in A_0$.
We need a technical lemma that follows from the proof of lemma 2.2.
In the following, for any point $w$\ in an algebraic set $W$, we will
use $\germ(W,w)$\ to denote the $\eps$-ball $B_{\eps}(w)$\ for
$0<\eps\ll1$\ of $w\in W$, under some Riemannian metric.

\pro{Lemma 2.4}
Let $\E_S$, where $s_0\in S$\ is a smooth curve, be
a family of $\mu$-unstable sheaves satisfying $\E_{s_0}\cong \bar\E$.
Then if $s_0\in U\sub S$\ and $\eta_1$\ and $\eta_2$\ are two analytic
maps from $(U,s_0)$\ to
$(A,0)$\ given by
(2) of the definition 2.1 based on the family $\E_S$, then the images
$\eta_1(\germ(U,s_0))$\ and $\eta_2(\germ(S,s_0))$\ are
contained in the same irreducible component of $\germ(A_0,0)$.
\endpro

\proof
We continue to use the notation developed in lemma 2.2.
Let $\e_U=\e_{S|C\times U}$.
Since $A\sub\Qu$, the map $\eta_i\mh U\to A$\
is given by quotient sheaf homomorphism
$f_i\mh p_Z\sta\R\to \E_{U}$. Since $f_i$\
is uniquely determined by the induced isomorphism $\tilde f_i\mh
\oplus^N\OO_{Z\times U}\to p_{U\ast}(\E_{U}\otimes p_Z\sta H)$\
(at least near $0\in\Qu$),
$\tilde g=\tilde f_2\upmo\circ \tilde f_1\in
\Hom(\OO^{\oplus N}_{Z\times U}, \OO^{\oplus N}_{Z\times U})$\
induces a homomorphism
$g\mh p_Z\sta\R\to p_Z\sta\R$\ that makes $f_1=f_2\circ  g$.
Since $\eta_1(s_0)=\eta_2(s_0)$\ as quotient sheaves,
$\tilde g(s_0)=c\cdot id$.
Next, let $T$\ be a connected analytic neighborhood of $1\in GL(N)$\ and
let $\Psi\mh T\times U\to \Qu$\ be the map defined by
$\Psi(h,s)=h\cdot f_2(s)$. By shrinking
$T$\ and $U$\ if necessary, we can assume
the composition of $\Psi$\ with
the local projection $\pi_{+}$\ of (2.2) is well-defined
and such that
$$\eta_i\Bl \germ(U,s_0)\Br \sub \pi_{+}\circ\Psi\Bl \germ(T\times U,
(1,s_0))\Br \sub \germ(A_0,0)\ \quad \text{for } i=1,2.
$$
Because $T$\ is irreducible, $\eta_1(\germ(U,s_0))$\ and
$\eta_2(\germ(U,s_0))$\ must be contained in
the same irreducible component of $\germ(A_0,0)$. This
completes the proof of the lemma.
\endpf

In the remainder of this section, we shall prove two technical
lemmas concerning $A_0$.

\pro{Proposition  2.5}
Assume $\eb$\ is locally free, then
the subvariety $A_0\sub A$\ has pure complex codimension $g+1$\ at
$0\in A_0$.
\endpro

When $\e$\ is a sheaf with torsion, we denote by $\e^t$\ the
torsion subsheaf of $\e$.

\pro{Proposition  2.6}
Assume $\ell(\eb^t)\leq 1$, then
the subvariety $A_0$\ is locally irreducible at $0\in A_0$.
\endpro

Here we say a variety $Z$\ is locally irreducible at $z\in Z$\
if for $0<\eps\ll 1$,
the smooth locus of the $\eps$-neighborhood $B\lep(z)\cap Z$\
is connected.
It is easy to see that if $w\in A_0$\ corresponds to a
locally free sheaf in $\s_1^0$, then $w$\ is a smooth point of $A_0$.
Also, an easy tangent space calculation shows that $\codim(\s^0_1\cap
A,A)=g+1$. Therefore, Proposition 2.5 will be established if we can
show that
$A_0\cap \s_1^0$\ is dense in $A_0$. To achieve this, we need to study
the deformation of sheaves in $\s_l$.

We first study the locally
free case. Let $\E\in\Sl$, $l\geq 2$, be locally free
and let $\L\sub\E$\ be the destabilizing subsheaf
of $\E$. We choose an effective divisor $D\sub C$\
of degree $l-1$\ and choose a homomorphism $\e\to\OO_D$\
so that the composition $\l\to\e\to\OO_D$\ is surjective.
Let $\f$\ be the kernel of $\e\to\OO_D$. Then
$\f$\ belongs to the exact sequence
$$0\lra \L(-D)\lra \F\lra \L^{-1}\lra 0
$$
and $\e$\ can be recovered as the kernel of
$\F(D)\mapright{\sigma_0}\OO_D$.
Next we fix a homomorphism
$\eta\mh \F(D)\to \OO_D$\ so that $\ker\{\eta\oplus \sigma_0\mh
\F(D)\to\OO_D^{\oplus 2}\}$\ is isomorphic to $\F$. Then for any
$s\in H^0(\OO_D)$, $\sigma_s=\sigma_0+s\eta: \F(D)\lra \OO_D$\
is always surjective and consequently, $\E_s=\ker\{\sigma_s\}$\
is locally free. $\{\E_s\mid s\in H^0(\OO_D)\}$\ provides us an algebraic
family of locally free sheaves parameterized by the total space $S$\ of
$H^0(\OO_D)$. We denote this family by $\E_S$.
In case $s\in S$\ corresponds to $0\in H^0(\OO_D)$,
then $\E_{s}$\ is isomorphic to $\E$. Further,
for any $s\in S$, $\E_s$\ is $\mu$-unstable because $\L(-D)\sub \E_s$\
and $\deg\L(-D)=1$. It is straight forward to check that for general
$s\in S$, $\l(-D)$\ is the destabilizing subsheaf of $\e_s$.
Therefore, we have proved

\pro{Lemma 2.7}
Let $D\sub C$\ be any effective divisor of degree $l-1$\ and let $\E_S$\
be the family of locally free sheaves  constructed with $S$\ the total
space of $H^0(\OO_D)$. Then $S\cap\s_1$\ is dense in $S$.
\endpro

Our next task is to analyze the set $S\cap \s_1^0$.
Let $\f$\ be given by
$$0\lra \L\lra \f\lra \L^{-1}\lra 0\tag 2.4
$$
with $\L$\ locally free and has degree $>0$.
We claim that when (2.4) does not split, then
$h^0(\endo(\f))= h^0(\L\uptwo)$.
Indeed, since $\deg\l>0$, there is a surjective $\Hom(\f,\f)\to
\Hom(\l\upmo,\l\upmo)=\CC$. Let $\Lam$\ be its kernel. Since
(2.4) is non-split, all elements in $\Lam$\ lift to $\Hom(\f,\l)$,
which lift to $\Hom(\l\upmo,\l)$\ for the same reason. Therefore,
$h^0(\Cal End(\f))=h^0(\l^{\otimes 2})+1$. This proves the claim.
Now let $s\in S\cap \s_1$\ be any point, where $S$\ is as in
lemma 2.7, and let $\E_s$\ be the locally free sheaf.
Then $\e_s$\ is non-split and then
$h^0(\endo(\E_s))=0$\  if and only if $h^0(\L\uptwo(-2D))=0$. To this end,
we need

\pro{Lemma 2.8}
For any $\L\in\Pic(C)$\ of degree $l\geq 2$, there is an effective divisor
$D\sub C$\ of degree $l-1$\ such that $h^0(\L\uptwo(-2D))=0$.
\endpro

\proof
Recall that for any subspace $V\sub H^0(\L\uptwo)$, there is an
$x\in C$\ so that
$\dim\bl V\cap H^0(\L\uptwo(-2x))\br \leq \dim V-2$.
Hence, the lemma follows if
$h^0(\L\uptwo)\leq 2\deg \L-2$. But this follows from the Clifford's
theorem and Riemann-Roch theorem.
\endpf

Combining lemma 2.7, 2.8 and the fact that $\e_s$\ is non-split
for general $s\in S$, we have proved

\pro{Proposition 2.9}
Let $\E\in\s_l$, $l\geq 2$\ be any locally free sheaf (with
$\det \E\cong \OO_C$). Then for
general effective divisor $D\sub C$\ of degree $l-1$\ and for the
family of locally free sheaves $\E_S$\ parameterized by $S$($=$\
the total space of $H^0(\OO_D)$) constructed in lemma 2.7, we have
$S\cap\s_1^0$\ is dense in $S$.
\endpro

\pro{Remark 2.10}
Let $\e\in\s_l$, $l\geq 2$\ be locally free,
let $Z$\ be the set of degree $l-1$\ effective
divisors (which is the symmetric product $S^{l-1}C$),  let
$D_Z\sub C\times Z$\ be the universal divisor and let $\tilde Z$\ be the total
space of the underlining vector bundle of the locally free sheaf
$\pi_{Z\ast}(\OO_{D_Z})$. Then all $\E_s$'s constructed in lemma 2.7
are parameterized by
$\tilde Z$. Further, if we let $Y\sub \tilde Z$\ be the subset corresponding
to the zero section of $\pi_{Z\ast}(\OO_{D_Z})$, then sheaves associated to
$w\in Y$\ are isomorphic to $\E$. Now let $\tilZ_0\sub\tilZ$\ be the set
of points whose associated sheaves belongs to $\s_1^0$. Then
lemma 2.9 says that $\tilZ_0$\ is dense in $\tilZ$.
\endpro

Now we study the deformation of locally free sheaves in $\s_1$. Let
$\E\in\s_1$\ and let $\L$\ be the destabilizing subsheaf of $\E$. We
choose an affine neighborhood $U$\ of $\{\l\}\in\Pic(C)$\ and let
$\L_U$\ be the restriction of the Poincare line bundle to $C\times U$.
We denote by $s_0\in U$\ the point $\l$. Now we construct a set that
parameterize sheaves having $\l_u$, $u\in U$, as their destabilizing
subsheaves. This can be done as follows:
Consider the extension group
$\Ext^1_{C\times U}(\L_U^{-1},\L_U)$. Since $U$\ is affine, this
group is an $\OO_U$-module and because $\Ext^2(\l\upmo_U,\l_U)=0$,
by cohomology and base change
theorem, the restriction homomorphism
$$ \Ext^1_{C\times U}(\L_U^{-1},\L_U)\lra \Ext^1_C(\L_u^{-1},\L_u)
\tag 2.5
$$
is surjective for each $u\in U$.
Let $Z$\ be the total space of $\Ext^1_{C\times U}(\L_U^{-1},\L_U)$.
There is an extension sheaf over $C\times Z$
$$0\lra p_{C\times U}\sta\l_U\lra \e_Z\lra p_{C\times U}\sta\l\upmo_U
\lra 0\tag 2.6
$$
such that whose restriction to each $C\times\{z\}$, $z\in Z$\ over
$u\in U$, is the extension sheaf defined by
$\tilde z\in \Ext^1_C(\L_u^{-1},\L_u)$\ that is
the image of $z$\ under (2.5).
Here $p_{C\times U}$\ is the projection $C\times Z\to C\times U$.
Because (2.5) is surjective, any extension sheaf of $\L_u^{-1}$\
by $\l_u$\ appears in this family. Therefore,
$\E_Z$\ is a deformation of $\E$\ whose destabilizing subsheaves are
members of $\L_U$. Of course, for general $z\in Z$, the sheaf $\e_z$\
is non-split.
Hence, $h^0(\endo(\E_z))=h^0(\l^{\otimes 2})=0$\ for general $z\in Z$\
because $g(C)\geq 4$.

With the material prepared, now it is easy to prove the
following proposition.

\pro{Proposition  2.11}
Let $\bar\E$\ be any $\mu$-unstable rank two locally free sheaf of
$\det \E\cong \OO_C$\ and let $A$\ be the versal deformation
space of $\bar\E$\ with $A_0\sub A$\ the locus of $\mu$-unstable sheaves. Then
$A_0\cap\s_1^0$\ is dense in $A_0$.
\endpro

\proof
It suffices to show that any $\mu$-unstable locally free sheaf $\e$\ admits
a deformation whose general member belongs to $\s_1^0$.
Such deformation was constructed in proposition 2.8 and (2.6).
\endpf

\noindent
{\sl Proof of proposition 2.5}.
Since $A_0\cap\s_1^0$\ is dense in $A_0$, it suffices to show that
for $w\in A_0\cap\s_1^0$\ close to $0$, $\codim(A_0\cap\s_1^0\ \text{at}\ w)
=g+1$. Let $\F_A$\ be the family associated to the versal deformation
space $A$\ of $\bar\E$, let $w\in A_0\cap \s_1^0$\ and let
$(A\pri,\F_{A\pri};0\pri,\F_w)$\ be the versal deformation space of $\F_w$.
Then by lemma 2.2,
there is an analytic neighborhood $U$\ of $w\in A$\ and
an analytic map $f\mh(U,w)\to(A\pri,0\pri)$\ such that
$(1_C\times f)\sta
 \F_{A\pri}$\ is (analytically) isomorphic to $\F_A$\ restricting to
$C\times U$.
Note, $U\cap\s_1^0=f^{-1}(A\pri\cap\s^0_1)$. Therefore proposition
2.5 follows from (a) $\codim (A\pri\cap\s_1^0,A\pri)=g+1$\ and (b)
$f$\ is a submersion at $w$.
We first prove (b). Since both $A$\ and $A\pri$\ are smooth,
it suffices to
show that the differential $df\mh T_wA\to T_{0\pri}A\pri$\ is surjective.
But this follows from the remark 2.3 since $w\in A_0$\ is close to 0.
(a) is a straight forward dimension counting argument based on the
exact sequence (2.4).
This completes the proof of proposition 2.5.
\endpf

We now turn our attention to the proof of proposition 2.6
for locally free $\eb$.

\noindent
{\sl Proof of Proposition 2.6}. We first consider the case where
$\bar\E\in\s_l$\ is locally free and $l\geq 2$. Following the remark 2.10,
there is a pair of irreducible varieties $Y\sub \tilZ$\
and a family of locally free sheaves $\E_{\tilZ}$\
on $C\times\tilZ$\ so that any deformation $\E_s$\ of $\bar\E$\ constructed
in lemma 2.7 is part of this family and for $w\in Y\sub\tilZ$,
$\E_w\cong \bar\E$. By
lemma 2.2, there is an analytic map
$f_w\mh \germ(\tilZ,w)\to (A,0)$\ associated to $w\in Y$\ that is
induced by the family $\E_{\tilZ}$. Since $\tilZ$\ is irreducible,
$f_w$\ defines a unique irreducible component of $\germ(A_0,0)$\ and
by lemma 2.4, this irreducible component
is independent of the choice of $f_w$\ and $w\in
Y$. We denote this component by $B_0$.
In the following, we will show that $B_0$\ is the
only irreducible component of $\germ(A_0,0)$.

Let $B$\ be any irreducible component of $\germ(A_0,0)$.
Since $B\cap \s_1^0$\ is dense in $B$, there is a smooth analytic curve
$s_0\in S$, an analytic map $\varphi\mh (S,s_0)\to (B,0)$\ so that
$\varphi(S-s_0)\sub B\cap \s_1^0$.
Let $\E_S$\ be the pull-back of $\F_A$. Because
$\varphi(S-s_0)\sub B\cap \s_1^0$\ and because destabilizing subsheaves of
unstable (locally free) sheaf is unique, there is an
invertible sheaf $\L_S$\ on $C\times S$\
so that for $s\in S-s_0$, $\L_s$\ is the destabilizing
subsheaf of $\E_s$. Hence there is a homomorphism $\L_S\to \E_S$\
that induces an exact sequence
$$0\lra\L_S\lra \E_S\mapright{\beta} \L_S^{-1}\otimes \JJ_{\Sigma}\lra 0,
\tag 2.7
$$
where $\JJ_{\Sigma}$\ is the ideal sheaf of a zero-scheme
$\Sigma\sub C\times S$\ supported on $C\times \{s_0\}$.
We claim that there is an effective divisor $D_S\sub C\times S$\
containing $\Sigma$\
flat over $S$\ such that $D_S$\ has degree $l-1$\ along fibers of
$C\times S\to S$. Indeed, for any $z\in \text{supp}(\Sigma)$,
since $\E_s$\ is locally free, $\JJ_{\Sigma,z}$\
is generated by, say  $f^z_1,f^z_2\in \OO_{C\times S,z}$.
Without loss of generality, we can assume

\noindent
(2.8) $f^z_1|_{C\times \{s_0\}}$\ generates the locally free part of
$\JJ_{\Sigma}|_{C\times\{s_0\}}$\ at $z$.

\noindent
Then the union of $\{f^z_1=0\}\sub
C\times S$\ for $z\in\text{supp}(\Sigma)$\ form a divisor $D_S$\ flat
over $S$\ near the fiber $C\times\{s_0\}$.
By shrinking $S$\ if necessary, we can assume
$D_S$\ is flat over $S$. We now check that it has degree $l-1$\ along
the fiber $C\times \{s_0\}$. Because of (2.8), by restricting
(2.7) to $C\times \{s_0\}$, we get
$$0\lra \L_{s_0}\lra \E_{s_0}\mapright{\beta_{s_0}}
\L_{s_0}^{-1}(-D_{s_0})\oplus \tau\lra 0,
$$
where $D_{s_0}=D_S\cap (C\times\{s_0\})$.
{}From this we see $\deg D_{s_0}=l-1$\ since $\E_{s_0}=\bar\E$\ has
destabilizing subsheaf of
degree $l$\ and $\deg \L_{s_0}=1$.

Now by taking the preimage of  $\L_S^{-1}(-D_S)\sub \L_S^{-1}\otimes
\JJ_{\Sigma}$\ under $\beta$\ in (2.7),
we get a locally free sheaf $\F_S$\ that belongs to the exact
sequence
$$0\lra \L_S\lra \F_S\lra \L_S^{-1}(-D_S)\lra 0\tag 2.9
$$
and the exact sequence
$$0\lra \F_S\lra \E_S\lra \T\lra0,
$$
where $\T=\L_S^{-1}\otimes \JJ_{\Sigma}/\L_S^{-1}(-D_S)$. Since $\JJ_{\Sigma}$\
is generated by $f^z_1$\ and $f_2^z$\ and $\l\upmo_S(-D_S)$\ is generated
by $f_1^z$\ at $z\in \text{supp}(\Sigma)$,
$\T$\ is a rank one locally free $\OO_{D_S}$-modules. Therefore,
imitating the argument proceeding lemma 2.7, we can recover $\e_S$\ by
$$0\lra \E_S\lra \F_S(D_S)\mapright{\sigma_S}\OO_{D_S}\lra 0.
$$
Because of (2.9), $\F_s$\
has degree 1 destabilizing quotient sheaf for any $s\in S$.

Now we show that  $\varphi(\germ(S,s_0))\sub B_0$. Let
$0\in T$\ be a smooth analytic curve and let
$$\tilde \sigma_{S\times T}: \F_S(D_S)\otimes_{\OO_{C\times S}}
\OO_{C\times S\times T}\lra \OO_{D_S\times T}
$$
be a surjective homomorphism whose restriction to $C\times S\times0$\
is $\sigma_S$. Let $\e_{S\times T}$\ be the kernel of
$\tilde \sigma_{S\times T}$.
By definition 2.1, this family induces an analytic map
$\eta\mh \germ(S\times T,(s_0,0))\to (A,0)$. Thanks to lemma 2.4,
we know that $\eta(\germ(S,s_0)\times \{0\})\sub B$\ and
because the restriction to $C\times\{s_0\}\times T$\ of $\e_{S\times T}$\
belongs to the family constructed in the Remarks 2.10, $\eta(\{s_0\}\times
\germ(T,0))\sub B_0$. Further, by choosing $\tilde \sigma_{S\times T}$\
generic, we can assume for $t\in T-0$, $\E_{(s_0,t)}\in \s_1^0$.
Therefore, the general points of
$\eta(\germ(\{s_0\}\times T,(s_0,0))$\ are smooth points
of $A_0$. On the other hand, general points of
$\eta(\germ(S,s_0)\times\{0\})$\ are smooth points of $A_0$\ as well.
Hence $B_0=B$\ because $S\times T$\ is
irreducible. This completes the
proof of Proposition 2.6 for $l\geq 2$.

It remains to show that $\germ (A_0,0)$\ is irreducible when
$\bar\e\in \s_1$.
Let $\l$\ be the destabilizing subsheaf of $\bar\e$\ and let $\e_Z$\
be the family of sheaves on $C\times Z$\ constructed in (2.6).
Because the set $Z_0=\{z\in Z\mid \e_z\cong\bar \e\}$\ is connected and $Z$\
is irreducible, the data $Z_0\sub Z$\ defines a unique irreducible
component $B\sub\germ(A_0,0)$. On the other hand, because all deformations
of $\bar\e$\ can be realized as a sub-family in $Z$, by
lemma 2.4, $\germ(A_0,0)\sub B$. Therefore, $A_0$\ is locally irreducible at 0.
\endpf

In the remainder of this section, we will study the deformation of
$\mu$-unstable non-locally free sheaves.
First, we collect some information of
non-locally free sheaves in $\s_l$.

\pro{Lemma 2.12}
Let $\e\in\s_l$\ be any non-locally free sheaf, then there is a
deformation of $\e$\ so that whose generic member are locally free in
$\s_l$. Further, if $l\geq 2$, then we can find deformation of $\e$\
so that whose generic member are non-locally free in $\s_1$.
\endpro

\proof
Let $\etor\sub\e$\ be the torsion subsheaf and let $\ef=\e/\etor$.
Then $\e=\ef\oplus\etor$. Let $S$\ be the total
space of $\Ext^1(\etor,\ef)$. Then $S$\ defines a family of sheaves
that are extension of $\etor$\ by $\ef$. The desired families
can be chosen as subfamilies of $S$\ easily.
\qed

\pro{Proposition 2.13}
Let $\e\in\s_1$\ be any sheaf with length one torsion subsheaf and
let $(A,0,\f_A,\e)$\ be the versal deformation space given by lemma 2.2. Then
$A_0$\ is locally irreducible at 0, where $A_0$\ is the locus of
$\mu$-unstable sheaves.
\endpro

\proof
The proof is similar to that of proposition 2.6.
Let $B$\ be an irreducible component of $\germ (A_0,0)$.
Then by the previous lemma, the general sheaves in $B$\ are locally free
sheaves in $\s_1^0$. Hence we can find a smooth curve $s_0\in S$\
and a morphism $f\mh(S,s_0)\to (A,0)$\ so that
$f(S-s_0)\sub B\cap \s_1^0$\
and are locally free. Let $\f_S$\ on $C\times S$\ be the pull
back of $\f_A$\ and let $\l_S$\ be the family of (invertible)
destabilizing quotient sheaves of $\e_S$. Then there is a zero scheme $Z\sub
C\times S$\ such that $\e_S$\ is the extension of $\l_S$\ by
$\l_S\upmo\otimes\JJ_{Z\sub C\times S}$. Because $\ell((\e_{s_0})^t)=1$,
$Z\cap \bigl(C\times\{s_0\}\bigr)=\{z\}$\ (reduced) for some $z\in C$.
Therefore,
$$Z=\{t+h(t,s)=0,\ s^k=0\},\quad h(s,t)\in (s,t^2)\OO_{C\times S,
(z,s_0)}
$$
for some $k\geq 1$, where $t$\ and $s$\ are analytic coordinate
of $z\in C$\ and $s_0\in S$\ respectively. Since $Z$\ and
$\{t=0,s^k=0\}$\ are contained in the family
$$Z_u=\{t+u\cdot h(t,s)=0,s^k=0\},\quad u\in \CC,
$$
any family of sheaves $\e_S\pri$\ that restricts to $\e_{s_0}$\ on
$C\times\{s_0\}$\ and is an extension of $\l_S$\ by $\l_S\upmo\otimes
\JJ_{Z_u \sub C\times S}$, $u\in\CC$, will be in $B$\ thanks the lemma 2.4
(at least near $s_0$). Therefore, $B$\ depends only on the choice
of $k$\ and the family $\l_S$. But
because $\Pic(C)$\ is smooth, $A$\ does not depend on $\l_S$\ either.
Hence for families $\e_S$\ whose corresponding zero schemes $Z$\
have identical degree $k$\ (considered as cycles) after restricting
to fiber $C\times\{s_0\}$, they identify the same irreducible component
of $\germ (A_0,0)$. We denote this component by $A_{[k]}$. However,
if we take a base change of $S$\ branched at $s_0$\ of ramification index $m$,
the new family certainly is still contained in $A_{[k]}$. Hence
$A_{[k]}=A_{[mk]}=A_{[m]}$. This proves that $\germ(A_0,0)$\ is irreducible.
\qed

\head 3. Deformation of locally free sheaves
\endhead

In this section, we will construct the homology mentioned in the
introduction that will lead us to the proof of the vanishing
of (0.9).

We fix $(H_0,I)$, a $C\in |2n_0H_0|$\ as in (1.1), a precompact
neighborhood $\C\ssub \nsq^+$\ of $H_0\in\nsq^+$\ and the integer
$N$\ given in lemma 1.3 and 1.5. In this section, we will prove

\pro{Proposition 3.1}
Let $d\geq N$\ and let $H\in\c$\ be $(H_0,I,d)$-suitable.
Then the subset $\w\sub\mhid$\ defined in (3) of lemma 1.7
(i.e. $\w$\ consists of shaves in $\mhid$\
restricting to unstable sheaves on $C$) satisfies
$$H_i\bl\w,\w\cap\s\mhid\br=0,\quad i\leq 2.\tag 3.1
$$
\endpro

To prove the theorem, we need to show that every
homology cycle $(\sig^i,\partial \sig^i)\to (\w,\w\cap\s\mhid)$\ is
homologous to zero. We will construct such homology directly by constructing
deformation of sheaves.

We first study deformation of sheaves following [La, OG1].
Let $W$\ be any reduced (quasi-projective)
scheme and let $\E_W$\ be a family of rank two torsion free sheaves
on $X\times W$\ flat over $W$\ so that $\wedget\E_W\cong\OO_{X\times W}$.
We assume that $\e_W$\ is locally free along $C\times W$\ and
for any $w\in W$, the restriction of $\E_w$\ to
$C$\ is unstable and whose destabilizing quotient sheaf has
degree $-1$. Because the destabilizing quotient sheaf is
unique, there is an invertible quotient sheaf $\e_{W|C\times W}
\to\L_W$\ on $C\times W$\ whose restriction to $C\times\{w\}$,
$w\in W$, realizes $\l_w$\ as the destabilizing quotient sheaf of
$\e_{w|C}$. Let $\inl\mh C\times W\to X\times W$\
be the obvious immersion and let $\E_W\lra \inl_{\ast}\L_W $\
be the induced (surjective) homomorphism. We denote its kernel
by $\f_W$. Then $\F_W$, $\E_W$\ and $\L_W$\
are related by the following exact sequence
$$0\lra \E_W(-C\times W)\lra\F_W\mapright{\sigma_W} \inl_{\ast}\L_W\upmo
\lra 0.\tag 3.2
$$
Of course, $\E_W$\ is determined by the section
$$\sigma_W\in H^0(C\times W, \hom(\F_W,\L_W\upmo)).
\tag 3.3
$$
Next we will enlarge the family $\E_W$\ by varying the sheaf $\L_W$\
and the homomorphism $\sigma_W$.
Let $J=\Pic^0(C)$\ and let $\P$\ be the normalized
Poincare line bundle so that $\P_{|p_0\times J}\cong\OO_J$,
where $p_0\in C$\ is fixed. (For $\lam\in J$, we will use $\P_{\lam}$\
to denote the bundle $\P_{|C\times\{\lam\}}$.) For the moment
we assume for some large $n$, the sheaf $\hom(\F_{W},\L_{W}\upmo)$\
belongs to the exact sequence
$$0\lra \hom(\F_{W},\L_{W}\upmo)\lra\N_{W} \lra \T\lra 0,
\tag 3.4
$$
where $\N_{W}=\OO_{C\times W}^{\oplus 2}(np_0\times W)$\
and $\t$\ is a family of torsion sheaves on $C\times{W}$\
flat over $W$. Let $p_{C\times W}$, $p_\ctimej$\ and $p_{W\times J}$\
be projections from $\ctimew\times J$\ to $C\times W$, $C\times J$\
and $W\times J$\ respectively.
We consider the following direct image sheaves on $W\times J$:
$$\A\wtj=p_{W\times J\ast}(p_{C\times W}\sta\N_W\otimes p\ctj\sta\P);
$$
$$\B\wtj=p_{W\times J\ast}(p_{C\times W}\sta\T\otimes p\ctj\sta\P).
$$
Both $\A\wtj$\ and $\b\wtj$\ are locally free. We
let ${\bold A}$\ be the vector bundle associated to $\A_\wtimej$\
and let $\opmo$\ be the tautological line bundle of the projective
bundle $\PP({\bold A})$\ of ${\bold A}$. We use the convention that
$\opmo$\ is a subsheaf of $\pi\sta\A_\wtimej$, where
$\pi\mh\PP({\bold A})\to\wtimej$\ is the projection. Composed with the induced
homomorphism $\pi\sta\a\wtj\to \pi\sta\b\wtj$, we get
$$f: \opmo\lra \pi\sta\B\wtj. \tag 3.5
$$
Since both $\a$\ and $\b$\ are locally free, for any closed $(w,\lam)\in
\wtimej$\ and
$$ v\in\ \text{the total space of}\ \opmo\ \text{over}\ (w,\lam),
$$
$v$\ corresponds to a section (unique up to scalars) $\phi\in
H^0(\OO_C^{\oplus 2}(np_0)\otimes\P_{\lam})$.
Further, $f(v)=0$\ if and only if the image of $\phi$\ in
$H^0(\T_{|C\times\{w\}}\otimes\P_{\lam})$\ is trivial, and
by the exactness of (3.4) it occurs exactly when $\phi$\ can be lifted
to a homomorphism $\tilde \phi\in H^0(\F_w,\L_w\upmo\otimes\P_{\lam})$.
Hence, the set of all non-trivial homomorphisms $\Hom_C(\F_w,
\L_w\upmo\otimes\P_{\lam})$\ modulo scalars is parameterized by
the scheme $f\upmo(0)\sub\PP({\bold A})$.
We let $\tilde Z$\ be $f\upmo(0)$\ endowed with reduced scheme
structure.

\pro{Lemma 3.2}
Let $\c_{W\times J}=p_{W\times J\ast}
(p_{C\times W}\sta(\f_W\dual\otimes\l_W\dual)\otimes p\ctj\sta\P)$.
Then there is a homomorphism
$$\alpha_{\tilde Z}:\opmo_{|\tilde Z}\lra \pi\sta\c_{W\times J}\tag 3.6
$$
so that for any closed $z\in\tilde Z$\ over $(w,\lambda)\in W\times J$,
the image
$$\text{Im}(\alpha_z)\sub \Hom_C(\F_{w|C},\L\upmo_w\otimes
\P_{\lam})
$$
induced by the restriction homomorphism
$$\pi\sta\c_{W\times J}\otimes k(z)
\lra \Hom_C(\F_{w|C},\L\upmo_w\otimes\P_{\lam})
$$
is non-trivial. Further, there is a section $\rho\mh
W\to \tilde Z$\ of the projection $\pi_W\mh \tilde Z\to W$\
so that for any closed
$w\in W$\ and $z=\rho(w)$, $\text{Im}(\alp_z)$\ coincides with
$\sigma_W\otimes k(w)$\ of (3.3) up to scalars.
\endpro

\proof
$\pi\sta\a_{W\times J}$, $\pi\sta\b_{W\times J}$\ and $\pi\sta\c_{W\times J}$\
certainly belong to the exact sequence
$$\pi\sta\c_{W\times J}\lra \pi\sta\a_{W\times J}\mapright{\beta}
\pi\sta\b_{W\times J}.
$$
Because the composition of
$$\opmo_{|\tilde Z}\lra
\pi\sta\bl\a_{W\times J}\br
$$
with $\beta$\ is trivial, it lifts to a unique homomorphism
$\alp_{\tilde Z}$\ in (3.6). The non-triviality of $\text{Im}(\alp_z)$\ for
$z\in \tilde Z$\ is obvious.
It remains to prove that there is a section
$\rho\mh W\to \tilde Z$\ with the desired property.
Indeed, if we let $\lam_0\in J$\ be the trivial
line bundle and let $i\mh W\mapto{\cong} W\times\{\lam_0\}
\sub W\times J$\ be the immersion, then the
section $\sigma_W$\ in (3.3) provides a section $\tilde\sigma_W\in
H^0(C\times W,\N_W)$\
that induces a section $\rho\mh W\to
\PP({\bold A} )$\ whose image is contained in $\tilde Z$.
$\rho\mh W\to \tilde Z$\ is the desired section.
\qed

In the following, we will use the section $\alpha_{\tilde Z}$\ to
construct a family of torsion free sheaves on $X\times \tilde Z$.
Let $\pi_{C\times W}$, $\pi_{C\times J}$\ and
$\pi_{\tilde Z}$\ be projections
from $C\times\tilde Z$\ to $C\times W$, $C\times J$\ and $\tilde Z$\
respectively.
Then with $p_{X\times W}\mh X\times \tilde Z\to X\times W$,
the section $\alpha_{\tilde Z}$\ provides us a homomorphism
(on $X\times \tilde Z$)
$$\alpha_{\tilde Z}: p_{X\times W}\sta \F_W\lra
\inl_{\ast}\Bigl(\pi_{\tilde Z}\sta(\oppo_{|\tilde Z})
\otimes \pi_\ctimew\sta\L\upmo_W \otimes \pi_{C\times J}\sta\P\Bigr),
\tag 3.7
$$
$\inl\mh C\times \tilde Z\to X\times\tilde Z$.
We denote the right hand side of (3.7) by $\inl_{\ast}\L_Z$, where
$\L_{\tilde Z}$\ is an invertible sheaf on $C\times \tilde Z$.
For technical reason, we let $Z\sub\tilde Z$\ be the union of
those irreducible components $A\sub\tilde Z$\ such that at general
$z\in A$, the restriction to $X\times \{z\}$\ of $\alp_{\tilde Z}$\ is
surjective. Note that $\rho(W)$\ is contained in $Z$. We let $\alp_Z$\ and
$\l_Z$\ be restriction of $\alp_{\tilde Z}$\ and $\l_{\tilde Z}$\ to $X\times
Z$\ and $C\times Z$\ respectively.
We let $\Sigma\sub C\times Z$\ be the subscheme so that $\text{Im}(
\alpha_Z)=\inl_{\ast}\L_Z\otimes\JJ_{\Sigma}$\ ($\JJ_{\Sig}$\
is the ideal sheaf of $\Sig\sub C\times Z$) and let
$\tilde\E_Z(-C\times Z)$\ be the
kernel of $\alp_Z$. Then we have the following exact sequence (on $X\times Z$)
$$0\lra \tilde\E_Z(-C\times Z)\lra p_{X\times W}\sta\F_W
\mapright{\alp_Z}\inl_{\ast}
\L_Z\otimes \JJ_{\Sigma}\lra0.
\tag 3.8
$$
We need a technical lemma which says that $\tile_Z$\ is a
family of torsion free
sheaves flat over $Z$. Once we know that it is flat over $Z$,
then $\tilE_Z$\ is a family that contains $\e_Z$\ as a subfamily.

\pro{Lemma  3.3}
Let $\tilE_Z$\ be the sheaf in (3.8). Then
$\tilE_Z$\ is a family of torsion free sheaves on $X\times Z$\ flat over
$Z$. Further, if we let $Z_0=\{z\in Z\mid \alpha_z\ \hbox{is
surjective}\}$, then for closed $z\in Z$, $\tilde\E_z$\ is locally
free along $C$\ (resp. $X-C$) if and only if $z\in Z_0$\ (resp. $\e_w$\ is
locally free, where $z$\ lies over $w\in W$).
\endpro

\proof
When $z\in Z_0$, then $\alpha_Z$\ is surjective near $C\times\{z\}$
and the kernel $\E_Z$\
is locally free and flat there. Now we assume
$z\in Z-Z_0$. Because $p_{X\times W}\sta\f_W$\ is locally free
at $C\times\{z\}$,
$\sig$\ is locally defined by at most two equations. On the
other hand, $\sig$\ contains no fibers of $C\times Z\to Z$ and
$Z_0$\ is dense in $Z$\ by our selection of $Z\sub\tilde Z$.
Therefore, $\codim(\sig)\geq 2$. Thus
$\sig$\ is a local complete intersection scheme of
codimension 2. Next, because (3.8) is exact and
because $\inls\L_Z\otimes\JJ_{\Sigma}$\ is flat over $Z$,
$\tilde\E_Z$\ is flat over $Z$. As to the torsion freeness, it suffices
to show that
$${\Tor}(\JJ_{\Sig},\OO_{C\times\{z\}})=0.
\tag 3.9
$$
Now we prove (3.9). Let $p\in \Sigma$\
be any point over $z$. Since $\JJ_{\Sigma,p}$\ is generated by two
sections, say $f_1,f_2$, $\JJ_{\Sigma,p}$\
belongs to the exact sequence
$$0\lra \JJ_{D_1\subset C\times Z,p}\lra \JJ_{\Sigma,p}\lra
\JJ_{D_2\subset D_1,p}\lra 0,
$$
where $D_1$\ is the divisor $\{f_1=0\}$\ in $C\times Z$\ and $
D_2$\ is the divisor $\{f_2=0\}$\ in $D_1$, and
$\JJ_{D_1\subset C\times Z}$\ and $\JJ_{D_2\subset D_1}$\
are ideal sheaves of $D_1\sub C\times Z$\ and $D_2\sub D_1$\
respectively. Without loss of generality, we can assume $D_1$\
is a divisor flat over $Z$\ (because $\supp(\Sigma)$\ contains no fibers of
$C\times Z\to Z$). Therefore, $\JJ_{D_1\sub C\times Z,p}$\
and $\JJ_{D_2\sub D_1,p}$\
are locally free $\OO_{C\times Z,p}$-module and $\OO_{D_1,p}$-module
respectively, and hence have trivial $\Tor(\cdot,\OO_{C\times\{z\}})$.
Thus, (3.9) holds and $\tilE_Z$\ is a flat family of torsion free sheaves.

Finally, for any $z\in Z$,
the restriction of (3.8) to $X\times\{z\}$\ is still exact:
$$0\lra \E_z(-C)\lra \F_w\lra \Bl\inls\L_z\otimes \JJ_{\Sigma}\Br_{
|X\times\{z\}} \lra 0.
\tag 3.10
$$
In particular, when $p\in C$, $\E_{z,p}$\ is locally free
if and only if $\bl\JJ_{\Sigma}\br_{|C\times \{z\}}$\ is a locally free
$\OO_{C\times \{z\}}$-modules. Therefore, $\E_z$\ is locally free at
$p\in C$\ if and only if $(p,z)\not\in
\Sigma$\ and locally free at $p\not\in C$\ if and
only if $\f_w$\ is locally free at $p$. This completes the proof of the lemma.
\qed

We have the following lemma which states that the previous
construction is independent of the choices made.

\pro{Lemma 3.4}
The scheme $\tilde Z=f\upmo(0)$, $Z\sub\tilde Z$\
and the family $\tilE_Z$\ are
independent of the choice
of the inclusion $\hom(\F_W,\L_W\upmo)\to\N_W$.
\endpro

\proof
We only need to check the following:
Assume $\F_W\dual\otimes\L_W\upmo\to
\OO^{\oplus 2}_{C\times W}(mp_0\times W)$\ is another inclusion
with cokernel $\T\pri$\ making the diagram
$$\CD
0 @>>> \hom(\F,\L_W\upmo) @>>>
\OO_{C\times W}\uptwo(np_0\times W) @>>> \T @>>> 0\\
@. @|  @VVV  @VVV \\
0 @>>> \hom(\F,\L_W\upmo) @>>> \OO_{C\times W}
\uptwo(mp_0\times W) @>>> \T\pri @>>> 0\\
\endCD
$$
commutative and assume the vertical arrows are injective, then
if we let $Z\sub\tilde Z\sub\PP({\bold A})$\ and $Z\pri\sub\tilde Z\pri
\sub\PP({\bold A}
\pri)$\ be the corresponding subschemes,
there is an isomorphism between $Z\sub\tilde Z$\ and
$Z\pri\sub\tilde Z\pri$\ (canonically) and isomorphism between
$\tilE_Z$\ and $\tilE_{Z'}$\ (non-canonically).
This is obvious because ${\bold A}$\ is a
subbundle of ${\bold A}\pri$\ and under the inclusion $\PP({\bold A})
\sub\PP({\bold A}\pri)$, $\tilde Z\pri\sub\PP({\bold A})$. We
shall omit the details here.
\qed

\pro{Remark}
Indeed, more is true: Let $W_i$, $i=1,2$, be two varieties with family $\e_i$\
on $X\times W_i$\ and let $Z_i$\ be the corresponding schemes constructed.
Assume there is a morphism $\varphi\mh W_1\to W_2$\
that induces an isomorphism
$$\Bigl((1_X\times \varphi)\sta\e_2\Bigr)
_{|2C\times W_1}\cong \e_{1|2C\times W_1},
$$
then there is a canonical isomorphism $Z_1\cong W_1\times_f Z_2$\ that
respects $\rho_i\mh W_i\to Z_i$.
The proof of it is similar to that of lemma 3.4. The key observation is that
the construction of $Z$\ involves twice elementary transformations of
sheaves $\e_i$\ along $C$\
which depends only on their restrictions to $2C$.
\endpro

We now study a subset of $Z$\ of non-locally
free sheaves: We keep $p_0\in C$\ and define
$$Z^{p_0}=\{z\in Z\mid \tilde\E_z\ \hbox{is not locally free at }\ p_0\}.
$$
Clearly, $Z^{p_0}\sub Z$\ is the the vanishing locus
of the restriction to $\{p_0\}\times Z$\ of (3.7)
$$\alpha_{Z,p_0}: \pi\sta\f_{0}\lra \OO_{{\bold P}}(1)_{|Z}
\otimes \bl \pi\sta_\wtimec\L_{W}\br_{|\{p_0\}\times Z},\tag 3.11
$$
where $\pi\mh Z\to W$\ is the projection and $\f_0\cong \f_{W|\{p_0\}\times
W}$\ is viewed as a sheaf on $W$. Because the restriction
of $\bl \pi\sta_\wtimec\L_{W}\br_{|\{p_0\}\times Z}$\ to any fiber of $\pi$\
is trivial and the line bundle $\oppo$\ is $\pi$-relatively ample [La,p46],
$Z^{p_0}\sub Z$\ is $\pi$-relatively ample
and has codimension at most 2.
In the following, we will apply the generalized Lefschetz hyperplane
theorem to the pair $(Z^{p_0},Z)$. Let $Z_{\text{reg}}$\ be
the largest open subset such that the restriction of $\pi$\ to
$\zreg$\ is smooth [Ha, p271].

\pro{Lemma 3.5}
Assume $\pi$\ is smooth at general points of $\rho(W)\sub Z$. Then there
is a maximal dense open subset $W_0\sub W$\ of which the
following holds:
Let $Z_0=\pi\upmo(W_0)$. Then $\rho(W_0)\sub \zreg$\ and the induced
pairs $Z_0\to W_0$,
$\zreg\cap Z_0\to W_0$\ and $\zpo\cap Z_0\to W_0$\ are
topological fiber bundles. Further, there is a
Riemannian metric on $Z$\ so that for any $\del$-neighborhood
$Z^{p_0,\del}$\ of $Z^{p_0}\sub Z$, $0<\del\ll1$, and any $w\in W_0$,
$$H_i\bl \zreg\cap\pi\upmo(w), \zreg\cap Z^{p_0,\del}\cap \pi\upmo(w)\br=0
\tag 3.12
$$
for $i\leq 3$.
\endpro

\proof
First we prove the existence of such an open subset $W_0\sub W$.
Since $\pi\mh Z\to W$\ is projective, there are Whitney
stratification $\s_Z$\ of $Z$\ and $\s_W$\ of $W$\ by algebraic
subvarieties such that $\pi$\ is a stratified map (see [GM] for definition).
Without loss of generality, we can assume $\zreg$, $Z^{p_0}$\ and
$\rho(W)\sub Z$\ are union of strata of $\s_Z$.
Let $W_0$\ be union of top-dimensional strata in $\s_W$\ and let
$Z_0=\pi\upmo(W_0)$. $W_0\sub W$\ is open and dense.
Because $\pi$\ is smooth at general points of $\rho(W)\sub Z$,
$\rho(W_0)\sub\zreg$. Further,
by Thom's first isotopy lemma, $Z_0\to W_0$, $\zreg\cap Z_0\to
W_0$\ and $Z_0^{p_0}\to W_0$\ are topological fiber bundles.
As to the vanishing of (3.12), we notice that by (1.1),
$\dim\pi\upmo(w)\geq 10$\ for any $w\in W$\ [La, p45].
On the other hand, $Z_0^{p_0}\cap \pi\upmo(w)$\
is the vanishing locus of two sections of the ample line bundle
$\oppo_{|Z_{w}}$. By the theorem on
page 150 and remarks on page 152 of [GM], we can find a
Riemannian metric on $Z$\ that provides us the vanishing of (3.12).
Here, we have used the fact that for $w\in W_0$, $\zreg\cap\pi
\upmo(w)$\ is smooth.

Finally, by using the fact that if $W_0, W_0\pri\sub W$\ are two open
subsets satisfying the conclusion of the lemma, then
$W_0\cup W_0\pri$\ also satisfies the conclusion of the lemma. Thus, there
is a maximal open subset $W_0\sub W$\ of which the lemma holds.
\qed

Now we are ready to construct the homology between any
class in $H_i(\w)$\ with class in $H_i(\w\cap\s\mhid)$. Let $e=-2C\cdot C$.
In light of the previous
construction, we will take $W$\ to be the set

$$ W=\left\{\E\in\w
\Bigm\vert
\matrix\format\l\\
   \text{$\e$\ is $e$-stable, $\e_{|C}$\ is locally free and
the destabilizing}\\
   \text{quotient sheaf of $\e_{|C}$\ has degree
${1\over 2}I\cdot C-1$.}
\endmatrix\right\}
\tag 3.13
$$

To construct the corresponding $\pi\mh Z\to W$\ and sheaves
$\tile_Z$\ on $X\times Z$, we argue as follows:
We first take an (analytic or \etale)
open covering $\{W_\alp\}$\ of $W$\ so that over each $W\lalp$, there is
a universal family $\E\lalp$\ on $X\times W\lalp$\ that
belongs to the exact sequence (3.4) for some integer $n\lalp$.
Then we apply lemma 3.2 and 3.4 to pairs $(W\lalp,\e\lalp)$\ to get
$\pi\lalp\mh \tilde Z\lalp\to W\lalp$\
and sheaves $\tilE\lalp$\ on $X\times \tilde Z\lalp$.
By lemma 3.4, when $W_{\alpha\beta}=W\lalp\cap W_
{\beta}\ne\emptyset$,
$\pi\lalp\upmo(W_{\alp\beta})$\ and $\pi_{\beta}\upmo
(W_{\alpha\beta})$\ are canonically isomorphic and
$\tilE\lalp$\ and $\tilE_{\beta}$\ are isomorphic (non-canonically)
on the overlap as well. Hence we can patch $\{\tilde Z\lalp\}$\
(resp. $Z\lalp\sub\tilde Z\lalp$;
$\rho\lalp\mh W\lalp\to Z\lalp$; resp. $Z\lalp^{p_0}\sub Z$)
together to get a scheme $\pi\mh\tilde Z\to W$\
(resp. $Z\sub \tilde Z$; $\rho\mh W\to Z$; resp.
$Z^{p_0}\to Z$). Also, the collection $\{\tilE\lalp\}$\ provides
us a rational map $\eta\mh Z--\to\mhid$\
by sending $z\in Z$\ to $\tilE_z$\ when it is $H$-stable.
Finally, because sheaves in $W$\ are $e$-stable, sheaves in $Z$\
are necessarily $H$-stable [OG1,p597]. Hence $\eta\to\mhid$\
is well-defined everywhere. Note that
$\eta\circ\rho$\ is the identity map,
$\eta({Z}^{p_0})\sub\s\mhid$\ and for any
$w\in W\cap\s\mhid$,
$\eta\bl\pi\upmo(w)\br\sub\s\mhid$.

\pro{Proposition 3.6}
Let $H_0$, $H_0\in\c\ssub\nsq^+$, $C\sub X$\ and $N$\ be as before.
For any $d\geq N$\ and $(H_0,I,d)$-suitable $H\in \c$,
let $W$\ be the set (3.13) and let $\pi\mh Z\to W$, $Z^{p_0}\sub Z$\
and $\rho\mh W\to Z$\ be the sets constructed. Then $W$\
is dense in $\w$\ and has pure codimension $g+1$\ in $\mhid$\ ($g=g(C)$).
Further there is a maximal
open dense subset $W_0\sub W$\ such that
\roster
\item
Let $Z_0=\pi\upmo(W_0)$, $Z^{p_0}_0=Z_0\cap Z^{p_0}$\ and
$Z_{0,\text{reg}}=Z_0\cap Z_{\text{reg}}$, then $\rho(W_0)\sub
\zoreg$\ and
$\pi_{|Z_0}\mh Z_0\to W_0$, $Z^{p_0}_0\to W_0$\ and
$Z_{0,\text{reg}}\to W_0$\ are
topological fiber bundles;
\item
There is a Riemannian metric on $Z$\ such that
for any $\del$-neighborhood ${Z}^{p_0,\del}$\ of
$Z^{p_0}\sub Z$, $0<\del\ll1$, and any $w\in W_0$,
$$H_i\bl \zreg\cap\pi\upmo(w),
Z^{p_0,\del}\cap Z\lreg\cap\pi\upmo(w)\br=0
$$
for $i\leq 3$.
\endroster
\endpro

\proof
We first prove that $\w$\ has pure codimension $g+1$. Let
$p\in\w$\ be any point corresponding to $\e$\ and let $U\sub\mhid$\
be an analytic neighborhood of $p$\ so that there is a well-defined
map
$$\varphi: (U,p)\lra (A,0)
\tag 3.14
$$
provided by lemma 2.2, where $(A,0)$\ is the versal deformation
space of $\e_{|C}$. Since $\varphi\upmo(A_0)=U\cap\w$\ and
$\codim(A_0,A)=g+1$, the codimension of each component of $\w$\ is
at most $g+1$. Further, when $p\not\in \Loc$\ (cf. (1.3)),
$\varphi$\ is a
submersion at $p$\ and the codimension of $\w$\ at $p$\ is exact $g+1$.
On the other hand, by choosing $N$\ large, we can assume $\dim\Loc\leq
\dim\mhid-10g$\ (lemma 1.3). Thus $\w-\Loc$\ is dense in $\w$\ and
therefore $\w$\ has pure codimension $g+1$\ as claimed.
As to the inclusion $W\sub\w$, it is clear that $\w=W\cup R_1\cup R_2$,
where $R_1$\ consists of sheaves in $\w$\ that are not $e$-stable and
$R_2$\ consists of sheaves in $\w$\ that are not locally free
along $C$. By lemma 1.1, there is an integer $N$\ such that for
$d\geq N$, $\dim R_1<\dim\w-4$. Since $\dim \Loc<\dim\w$, $W$\
will be dense in $\w$\ if we can show that any sheaf $\e\in R_2-\Loc$\
admits a deformation whose generic member is in $W$ (i.e. whose
restriction to $C$\ is locally free and unstable). But this
follows from lemma 2.12 and [GL,p84-87].

We now prove that $\pi\mh Z\to W$\ is smooth at general points of
$\rho(W)$. Let $w\in W$\ be any point corresponding to $\e_w$\
and let $\f_w$\ be the sheaf in (3.2). Namely, $\f_w$\ is the
kernel of $\e_w\to \l_w$, where $\l_w$\ is the destabilizing quotient
sheaf of $\e_{w|C}$. If we assume $w\in W$\ is general, then
$\e_w$\ will be $2e$-stable. Hence $\f_w$\ will be
$e$-stable. In this way, we can define a morphism $f\mh W\to\MM$\
that sends $\e_w$\ to $\f_w$, where $\MM$\ is a moduli scheme of
stable sheaves with appropriate Chern classes. We let the image scheme
be $Y$. From the construction of $Z$, we know that a neighborhood
of $\rho(w)\in\pi\upmo(w)$\ is isomorphic to the set of surjective
homomorphisms $\alp\mh\f_w\to\l\pri$\ up to scalars,
where $\l$\ is an invertible sheaf on $C$\
with $\deg\l=\deg\l_w$. Let $\e_{\alp}$\ be
the kernel of $\alp$. We claim that if $\beta\mh\f_w\to\l''$\ is another
homomorphism and $\e_{\beta}$\ its kernel,
then $\e_{\alp}\cong\e_{\beta}$\ if and only if $\l'\cong
\l''$\ and $\alp=\lambda\beta$\ for some $\lam\in\CC$. Indeed, the
isomorphism $\e_{\alp}\cong\e_{\beta}$\ will induce homomorphism
$h\mh\f_w(-C)\to\f_w$\ whose determinant $\det h$\ vanishes along $C$.
Because $\f_w$\ is $e$-stable, $h$\ must be of trace type, namely
$h(x)$\ is a multiple of identity homomorphism at each $x\in X$. Thus
$h$\ must vanishes along $C$\ and thus we get the following
commutative diagram
$$\CD
0 @>>> \e_{\alp} @>>> \f_w @>\alp>> \l' @>>> 0\\
@. @V{\cong}VV @Vh(-C)VV @VVV\\
0 @>>> \e_{\beta} @>>> \f_w @>\beta>> \l'' @>>> 0.\\
\endCD
$$
Hence $h(-C)=\lam\cdot\text{id}$, $\l'\cong\l''$\ and $\alp=\lam\beta$\ as
claimed.

Now let $z\in Z$\ be a point over $w\in W$, corresponding to sheaf
$\e_z$\ and $\e_w$\ respectively. Then by sending $z$\ to
$(\e_z,\e_w)\in W\times_Y W$\ we obtain a morphism
$$\Cal H: Z\lra W\times_Y W.
$$
By previous argument, when $w\in W$\ is a general point and $z\in Z$\
is near $\rho(w)$, $\Cal H$\ is one-to-one. On the other hand, since all
schemes involved are with reduced scheme structures,
$\Cal H$\ must be an isomorphism at $\rho(w)$\ for general
$w\in W$. Now let
$\text{pr}_2\mh W\times_YW\to W$\ be the projection onto the
second copy. It is straight forward to check that at general $w\in W$,
both $\text{pr}_2$\ and $W\times_YW$\ are smooth at $(w,w)$.
Finally, because the $\pi\mh Z\to W$\ is exactly $\text{pr}_2\circ
\Cal H$\ (at least near $(w,w)$), $\pi$\ will be smooth at
$\rho(w)\in Z$\ as
well. This proves that $\pi$\ is smooth at general points of $\rho(W)$.

The last statement follows from lemma 3.5
\qed

In the following, we will show how to use the pair $(W_0,Z_0)$\ to get the
desired vanishing theorem. We begin with a simple version of theorem
3.1.

\pro{Proposition 3.7}
Let the notation be as before, then for $i\leq 2$,
$$H_i\bl W_0,W_0\cap\s\mhid\br\lra H_i\bl\w,\w\cap\s\mhid\br
\tag 3.15
$$
is trivial.
\endpro

\proof
Let $\xi\in H_i(W_0,W_0\cap\s\mhid)$\ be any element represented
by $\sig\sub W_0$\ with $\d\sig\sub W_0\cap\s\mhid$. Let $\rho\mh
W_0\to Z_0$\ be the section and let $\sig'=\rho(\sig)\sub \zoreg$.
Because for each $w\in W_0$, the pair
$$\bl\zoreg\cap Z_0^{p_0,\del}\cap\pi\upmo(w),
\zoreg\cap\pi\upmo(w)\br
$$
is 2-connected for $0<\del\ll 1$,
we can find an $(i+1)$-chain $T'\sub\zoreg$\ whose
boundary $\d T'$\ has a decomposition $\d T'=\sig'\cap
A_1'\cap A_2'$\ such that $A'_1\sub Z_0^{p_0,\del}$\ and
$$A'_2\sub \bigcup\{\pi\upmo(w)\mid w\in\d\sig'\}.
\tag 3.16
$$
Because
$Z_0\to W_0$\ is proper, we can find $T\sub Z_0$\ so that $\d T=
\sig\pri\cup A_1\cup A_2$\ with $A_1\sub Z_0^{p_0}$\ and
$A_2$\ satisfies (3.16) with $A_2\pri$\ replaced by $A_2$.
Thus $\d(\eta(T))=\sig$\ {\sl modulo}
$\s\mhid$. In particular, the image of $\xi$\ in $H_i(\w,\w\cap\s\mhid)$\
is trivial.
\qed

It remains to show that (3.15) is surjective. We first let $S_1=\Lambda_1^C$.
By proposition 3.6, $\w\sub\mhid$\ has pure codimension $g+1$\
and $S_1\sub\w$\ has codimension at least $6g$. Note that $\mhid$\ is a
local complete intersection and $\w\sub\mhid$\ is defined by $3g-2$\
equations. Hence by lemma 1.10,
$$H_i\bl\w_1,\w_1\cap\s\mhid\br\lra H_i\bl\w,\w\cap \s\mhid\br
$$
is surjective for $i\leq 2$, where $\w_1=\w-S_1$.

\pro{Lemma 3.8}
With the notation as before and let $S_2=(\w_1-W_0)\cap\s\mhid$,
then for $i\leq 2$\ the induced homomorphism
$$H_i\bl\w_1-S_2,(\w_1-S_2)\cap\s\mhid\br\lra H_i\bl\w_1,\w_1\cap\s\mhid\br
\tag 3.17
$$
is surjective.
\endpro

\proof
We first let $S_3$\ be those $\e\in\w_1$\ so that
$\e\in\s\mhid-\so\mhid$\ and $\e_{|C}$\ is not locally free.
$S_3\sub\w_1$\ is closed and has codimension 3.
By the proof of lemma 1.10, (3.17) is surjective with $S_2$\ replaced by $S_3$\
because $\w_1-\s\mhid$\ is locally irreducible.
Let $\w_2=\w_1-S_3$\ and let $S_4=(\w_2-W_0)\cap \s\mhid$.
We claim that $\w_2$\ is locally irreducible
and transversal to $\s\mhid$. Let $p\in \w_2$\ be a point
associated to $\bar\e\in\so\mhid$\ and let $0\in A$\
be a versal deformation space of $\bar\e_{|C}$\ given in lemma 2.2.
Then the map $\varphi$\ in (3.14) is a submersion at $p$\ and
hence $\w_1$\ is locally irreducible at $p$\ since
$A_0$\ is locally irreducible (proposition 2.6). Also
$\s\mhid$\ and $\w_1$\ are transversal at $p$\ as stratified
sets. This is because the tangent space of fiber of $\varphi$\ is the
image $\Ext^1(\e,\e(-C))^0\to\Ext^1(\e,\e)^0$\ and the tangent space
of $\so\mhid$\ is the image of $H^1(\e\dual\otimes\e)$, and they
together span $\Ext^1(\e,\e)^0$. This proves that $\so\mhid\cap \w_1$\
is locally irreducible.

Next, we claim that $(\w_2-S_4)\cap\s\mhid$\ is dense in $\w_2\cap\s\mhid$.
Indeed, because of the vanishing of $\Ext^2(\e,\e(-2C))^0$,
a general sheaf $\e\in\w_2\cap\s\mhid$\ is locally free along
$C\sub X$\ whose restriction to $2C$\ is isomorphic to $\f_{|2C}$\
for a general $\f\in\w$.
Then by the remark after lemma 3.4 and the maximal of
$W_0\sub W$, $\e$\ is contained in $(\w_2-S_4)\cap\s\mhid$.
This establishes the claim.
Because $\w_2$\ is locally irreducible, is transversal to $\s\mhid$\
and because $(\w_2-S_4)\cap\s\mhid$\ is dense in $\w_2\cap\s\mhid$, the
following homomorphism
$$H_i\bl\w_2-S_4, (\w_2-S_4)\cap\s\mhid\br
\lra H_i\bl\w_2,\w_2\cap\s\mhid\br
$$
must be surjective for $i\leq 2$. This proves the lemma because
$\w_2-S_4=\w_1-S_2$.
\qed

Now we prove theorem 3.1.
We only need to show that (3.17) is surjective.
Let $\w_3=\w-\Lambda_1^C-(\w-W_0)\cap \s\mhid$. We already know that
$$H_i(\w_3,\w_3\cap\s\mhid)\to H_i(\w,\w\cap\s\mhid)
$$
is surjective. Because $\w_3\cap\s\mhid=W_0\cap\s\mhid$\ and
because $\w_3$\ is locally irreducible, the case for $H_1$\ follows
from lemma 1.11. Now let $\xi\in H_2(\w_3,\w_3\cap\s\mhid)$.
By lemma 1.11 again, we can find a representative $\sig\sub\w_3$\ with
$\d\sig\sub\w_3\cap\s\mhid)$\ such that $\sig\cap(\w_3-W_0)$\ is
discrete. Since the codimension of $\Lambda_2^C$\ in $\mhid$\
is at least $3g$\ (lemma 1.3), by
perturb $\sig$\ as we did in lemma 1.10, we can assume
$\sig\cap(\w_3-W_0)\cap\Lambda_2^C=\emptyset$.
In the following, we will construct a new representative of
$\xi$\ that is contained in $W_0$.

We still keep the representative $\Sig$. Let
$p\in \sig\cap(\w_3-W_0)$\ be associated to $\e$. Since $p\not\in
\Lambda_2^C$, the set
$\Lambda_{\e}^C\sub\mhid$\ (see (1.5)) has pure
codimension at most $3(-\xx(\OO_{2C}))+(g+1)$. Let $R$\ be an
irreducible component of $\Lambda_{\e}^C$\ containing $p$.
Then by a recent work of [OG2], we can choose $N$\ (depending only
on $C$) such that whenever $d\geq N$, then
$\overline{R}\cap\s\mhid\ne\emptyset$.
On the other hand, by [GL,p80] we indeed have $R\cap\so\mhid\ne\emptyset$.
Note that by dimension count, $R\cap\so\mhid-\Ltc\ne\emptyset$.
Now we pick a differentiable path $\rho\mh [0,1]\to R$\ connecting $p$\ and
$\tilde p\in R\cap\so\mhid-\Ltc$\ and let $\u$\ be a (classical)
neighborhood of $\rho([0,1])$\ in $\mhid$.
Without loss of generality, we can assume there is a universal family
$\e_{\u}$\ on $X\times \u$. Now let $0\in A$\ be the versal deformation
space of $\e_{|2C}$. Based on the proof of lemma 2.2, we can find an analytic
map $\varphi\mh (\u,\rho([0,1]))\to (A,0)$\ induced by the family $\e_{\u}$,
after shrinking $\u$\ if necessary. By further shrink $\u$\
if necessary, we can assume $\varphi$\ is a submersion at $\rho([0,1])$\ and
thus realize $\u$\ as a product $(-\eps,1+\eps)\times\u_0$\ with $\varphi$\
factor through $\varphi_0\mh \u_0\to A$\ ($\{t\}\times\u_0$\ is a normal
slice of $\rho([0,1])$\ at $\rho(t)$). Thus by the
remark after lemma 3.4 and the maximal of $W_0\sub W$,
we can assume without loss of generality that
$W_0\cap \u$\ is a (topological) fiber bundle
over $\varphi([0,1])\sub A$. Finally, because $\so\mhid$\ is transversal to
$R$\
at $\tilde p$, we can apply the technique in the proof of lemma 1.11
to find a 2-chain $T$\ contained in $W_0\cap\u$\ with $\d T=\Gamma_1\cup
 \Gamma_2$\ such that $\Gamma_2\sub W_0\cap \so\mhid$\ and $\Gamma_1\sub
\sig$\ is the boundary of $B\leps(p)\cap \sig$.
(Indeed, $T$\ can be made to be homeomorphic to $\Gamma_1\times[0,1]$\
and is the result of pulling $\Gamma_1$\ to $\u\cap\s\mhid$\ along
$W_0\cap\u$.) Hence $\Sig\pri=(\Sig-B\leps(p))\cup T$\
represents the homological cycle $[\sig]$\ while
$\sig\pri\cap(\w_3-W_0)$\ is one point less than $\sig\cap (\w_3-W_0)$.
Therefore, by iterating this process for each point in $\Sig
\cap(\w_3-W_0)$, we
eventually get a $\tilde\sig\sub W_0$\ that represents $[\sig]\in H_2(\w,\w\cap
\s\mhid)$. This shows that (3.15) is surjective. This and proposition 3.7
together prove the proposition 3.1.
\qed

\head 4. Proof of the main theorems
\endhead

In this section, we will first prove theorem 0.4 by using
Lefschetz hyperplane theorem. After that, we will study the pair
$\s\mdh\sub\mdh$\ in detail to establish both theorem
0.5 and 0.1.

We will keep the notation developed in the previous sections.
We fix $H_0$, $H_0\in\c\ssub\nsq^+$\ and the $N$\ given in lemma 1.3 and
1.5. For $d\geq N$, we pick a $(H_0,I,d)$-suitable $H\in\c$.
$\M_{H_0}(I,d\pri)$\ is birational to $\M_H(I,d\pri)$,
$d\pri=d$\ or $d+1$, which induces the following commutative square
$$\CD
H_i(\M_{H_0}(I,d)^0) @>\cong>> H_i(\mhid^0)\\
@V \tau(d)_i VV  @V \tilde\tau(d)_i VV\\
H_i(\M_{H_0}(I,d+1)^0) @> \cong >> H_i(\M_H(I,d+1)^0)\\
\endCD
$$
for $i\leq 2$. Following the discussion in the introduction, theorem 0.1
follows from the surjectivity of $\tilde\tau(d)_i$, which is
equivalent to the following theorem.

\pro{Theorem 4.1} Let $d\geq N$\ and $H\in\c$\ be $(H_0,I,d)$-suitable.
Then the homomorphism
$\tilde \tau(d)_i\mh
H_i(\mhid^0)\to H_i(\M_{H}(I,d+1)^0)$\ is surjective.
\endpro

As explained in the introduction, theorem 4.1 will be proved in two steps:
The first is to use the Lefschetz hyperplane theorem to prove the
following:

\pro{Theorem 4.2}
With the choice of $d$\ and $H$\ made in theorem 4.1, then for $i\leq 2$,
$$H_i\bl\mhid,\s\mhid\br =0.
$$
\endpro

\proof
Let $\Lam=\mbhid-\mhid$. Since $H$\ is $(I,d)$-generic,
$\Lam\sub\s\mbhid$. Also, since $\codim(\Lam)\geq 10g$\ (lemma 1.7),
and $\mhid$\ is a local complete intersection, by lemma 1.10,
$$H_i\bl\mhid,\s\mhid\br\lra H_i\bl\mbhid,\s\mbhid\br
$$
is an isomorphism for $i\leq 2$. For the same reason, the subset $\w\sub\mhid$\
(defined in (1.1)) and its closure $\overline{\w}$\ in $\mbhid$\ induces a
surjective homomorphism
$$H_i\bl\w,\w\cap\s\mhid\br\lra H_i\bl\wb,\wb\cap\s\mbhid\br
$$
for $i\leq 2$. Since $H_i(\w,\w\cap\s\mhid)=0$\ (proposition 3.1),
theorem 4.2 follows from the surjectivity of the homomorphism
$$H_i\bl\wb,\wb\cap\s\mbhid\br\lra H_i\bl\mbhid,\s\mbhid\br
\tag 4.1
$$
for $i\leq 2$.

For simplicity, in the following we will denote
by $\mm$\ the set of all $\mu$-stable
sheaves $\e$\ with $\Ext^2(\e\dual,\e\dual)^0=0$\ and by $\mmb$\ the
space $\mbhid$. Let
$$\Psi:\mmb\lra \PP^R\tag 4.2
$$
be the morphism constructed in lemma 1.5 and let $V\sub\PP^R$\
be the codimension $3g-2$\
linear subspace such that $\Psi\upmo(V)\cap\mm=\w\cap\mm$.
For $\del>0$, we let $\Vdel\sub\PP^R$\ be the $\del$-neighborhood
of $V\sub \PP^R$\ under the Fubini-Study metric and
let $\wbdel=\Psi\upmo(\Vdel)$\ and $\wdel=\Psi\upmo(\Vdel)\cap\mm$.
We first consider the triple
$(\mm,\s\mm,\s\mm\cap\wdel)$\ and its induced long exact sequence
$$H_i(\mm,\s\mm\cap\wdel)\mapright{\alp} H_i(\mm,\s\mm)
\mapright{\beta} H_{i-1}(\s\mm,\s\mm\cap\wdel).
\tag 4.3
$$
Let $U\sub\mm$\ be the open subset consisting of sheaves $\e\in\mm$\
such that $\ell(\e\doubledual/\e)\leq 4$. $U$\ is smooth, the
compliment of $U$\ in $\mmb$\ has codimension 5 and for any $u\in U$,
$\dim \Psi\upmo\bl\Psi(u)\br\leq 12$. Then by lemma 1.10,
$$H_i(U,U\cap\s\mm\cap\wdel)\lra H_i(\mm,\s\mm\cap\wdel)\tag 4.4
$$
is surjective for $i\leq 2$\ because $\s\mm\sub\mm$\ is
Cartier. Next, because the fibers of
$\Psi_{|U}\mh U\to\PP^R$\ have dimension at most 12 and $U$\ has
pure dimension much bigger than $3g+20$, we can
apply the stratified Morse theory technique (exactly the same as in the proof
of theorem 4.1 on page 195 of [GM]) to the map
$\Psi_{|U}$\ to conclude that
$$H_i(U\cap\wdel,U\cap\s\mm\cap\wdel)\lra H_i(U,U\cap\s\mm\cap\wdel)
$$
is surjective for $i\leq 2$. Then by (4.3) and (4.4), for $i\leq 2$
$$H_i(U\cap\wdel,U\cap\s\mm\cap\wdel)\lra H_i(\mm,\s\mm)
\tag 4.5
$$
will be surjective if $H_{i-1}(\s\mm,\s\mm\cap\wdel)=0$.

We claim that for $j\leq 1$,
$H_j(\s\mm,\soo\mm)=0$. $H_0=0$\ because $\soo\mm$\ is dense
in $\s\mm$\ (lemma 1.1). For $H_1$, let $f\mh([0,1],\partial[0,1])\to
(\s\mm,\soo\mm)$\ be a continuous map.
Since $\soo\mm$\ is dense in $\s\mm$\ (by lemma 1.11),
we can assume without loss of generality that $f\upmo(\s\mm-\soo\mm)$\ is a
finite set, say $\{p_1,\cdots,p_k\}$. Because $\s\mm$\ is a local complete
intersection
and the compliment of $\so\mm\cup\sto\m\sub\s\mm$\ has codimension
2 by lemma 1.10, we can choose $f$\ so
that the points $p_i$\ are all contained in $\st\mm$.
Let $R_1=\sto\mm$\ and let $R_2=\st\mm-\sto\mm$.
Since $R_1$\ is irreducible and dense in $\st\mm$,
by lemma 1.12 we can choose $f$\
so that all $p_i$\ belong to $R_2$. On the other hand,
because $\s\mm$\ is locally irreducible at $R_2$\ (this can
be checked directly), we
can perturb $f$\ within $\bigcup\bl B\ldel\bl f(p_i)\br
\cap\s\mm\br$\ to obtain a representative $f\pri$\ of $[f]$\ whose image
is contained in $\soo\mm$. Therefore, $[f]=0$\ and hence
$H_1(\s\mm,\soo\mm)=0$.

Now we consider the triple $(\s\mm,\soo\mm,\soo\mm\cap\wdel)$\
and its induced long exact sequence ($j\leq 1$)
$$\lra H_j(\soo\mm,\soo\mm\cap\wdel)\lra H_j(\s\mm,\soo\mm\cap\wdel)\lra
H_j(\s\mm,\soo\mm)=0.
$$
Because $H_j(\s\mm,\soo\mm\cap\wdel)\to H_j(\s\mm,\s\mm
\cap\wdel)$\ is surjective ($\soo\mm\sub\s\mm$\ is dense), by the above
exact sequence to show
$H_j(\s\mm,\s\mm\cap\wdel)=0$\ it suffices to show that $H_j(\soo\mm,\soo
\mm\cap\wdel)=0$. For this, we will look at the restriction to $\soo\mm
$\ of $\Psi$\ in (4.2),
$$\Psi\pri: \soo\mm\lra \PP^R.
$$
Because fibers of $\Psi\pri$\ has dimension 1, we can apply theorem on
page 153 of [GM] to conclude that $H_j(\soo\mm,\soo\mm\cap\wdel)=0$\
for $j\leq 1$\ and then $H_j(\s\mm,\s\mm\cap\wdel)=0$.
Therefore, we have proved the surjectivity of (4.5)
for $i\leq 2$\ which combined with lemma 1.10 yields the
surjectivity of
$$H_i(\wbdel,\s\mmb\cap\wbdel)\lra H_i(\mmb,\s\mmb),\quad i\leq 2.
\tag 4.6
$$
Finally, since both $\mmb$\ and $\s\mmb$\ are complete and (4.6)
holds for all $0<\del\ll1$, we obtain the surjectivity of (4.1)
by applying proposition 4.A.1 on page 206 of [GM]. This completes
the proof of the theorem.
\qed

In the remainder of this section, we will use theorem
4.2 to establish theorem 4.1. But first, we will fill in the details of the
definition of the homomorphism (0.5). For any $x\in X$, we let
$\sox\mm\sub\so\mm$\ be the set of $\e\in\so\mm$\ that are
non-locally free at $x$. Note that $\sox\mm$\ is a $\PP^1$-bundle
over $\mmmo$. Let $V_0$\ be a general fiber of this bundle.
Then the inclusion $V_0\sub\mm$\ and the bundle $\sox\mm\to\mmmo$\
induce the commutative diagram
$$\minCDarrowwidth{14pt}
\CD
0 @>>> H_i(V_0) @>>> H_i(\so^x\mm) @>>> H_i(\mmmo) @>>> 0 \\
@. @| @V r(d)_i VV @. @.\\
0 @>>> H_i(V_0) @>>> H_i(\mm) @<<< H_i(\mm^0),\\
\endCD
\tag 4.7
$$
where $i\leq 2$.

\pro{Lemma 4.3}
Let $d\geq N$.
For $i\leq 2$, the induced homomorphism $H_i(\mmo)\to H_i(\mm)$\
is injective and whose image $H_i(\mmo)\util$\ is a (linear)
compliment of the image $H_i(V_0)\util$\ of $H_i(V_0)\to
H_i(\mm)$. In particular, since the top row of (4.7)
is exact, we get a unique
homomorphism $\rho(d)_i\mh H_i(\mmmo)\to H_i(\mmo)$\ that coincides with
Taubes' homomorphism $\tilde\tau(d)_i$.
\endpro

\proof
We need to show that $H_i(\mm^0)\util$\ is a compliment of $H_i(V_0)\util$.
Let $i=1, 2$. We first show
that $H_i(\mm)\to H_i(\mm,\mmo)$\ is surjective. Since $\mm$\ is
smooth, $H_i(\mm)\cong H_i(\mmo\cup\so\mm)$. Thus it suffices to
show
$$H_i\bl\mmo\cup\so\mm\br\lra H_i\bl\mmo\cup\so\mm,\mmo\br
\tag 4.8
$$
is surjective. Let $\u$\ be a tubular neighborhood of $\so\mm
\sub\mmo\cup\so\mm$\ such that $\u-\so\mm$\ is an
$(\RR^4-0)/\ZZ_2$-bundle over $\mmmo$. Let $\xi\in H_i(\mmo\cup\so\mm,
\mmo)$\ be represented by $\sig\sub \bl\mmo\cup\so\mm\br\cap\u$.
Then $\d\sig$\ represents a trivial cycle in $H_{i-1}(\u)=
H_{i-1}(\u-\so\mm)$. Thus we can find another chain $\sig'\sub\u-\soo\mm$\
with $\d\sig'=-\d\sig$. Hence $\sig\pri\cup\sig$\ represents a
cycle $\tilde\xi\in H_i\bl\mmo\cup\so\mm\br$\ whose image in
$H_i\bl\mmo\cup\so\mm,\so\mm\br$\ is $\xi$. Thus (4.8) is surjective
and hence $H_i(\mm^0)\to H_i(\mm)$\ is injective.

On the other hand, the composition
$H_i(V_0)\to H_i(\mm)\to H_i(\mm,\mmo)$\
is an isomorphism because the intersection $(V_0,\s\mm)$\ in $\mm$\
is -2 (this is well defined because $V_0$\ is proper) and $\s\mm$\ is
irreducible. Therefore, $H_i(\mmo)\util$\ is a compliment of $H_i(V_0)\util$\
and the homomorphism $H_i(\mmmo)\to H_i(\mmo)$\ is well-defined.
It is straight forward to check that the homomorphism $\rho(d)_i$\
coincides with that of Taubes'.
\qed

\pro{Corollary 4.4}
Assume $d\geq N$, then the homomorphism $\tilde \tau(d)_i$\
is surjective if and only if the homomorphism
$$r(d)_i: H_i(\sox\mm)\lra H_i(\mm)\tag 4.9
$$
is surjective.
\endpro

Before we prove the surjectivity of (4.9),
let us state a technical lemma whose proof will be postponed
until the end of this section. For convenience,
in the following whenever a space $Z$\ admits an obvious map
$Z\to \mm$, then we will denote by $H_i(Z)\util$\ the image
of $H_i(Z)\to H_i(\mm)$.

\pro{Lemma 4.5}
Let $\f\in\s_j\mm$, $j\leq 3$, be any sheaf and let $V(\f)\sub\s_j\mm$\
be the set of those $\e$\ such that $\f\dual\cong\e\dual$\ and
$\ell\bl(\e\doubledual/\e)\otimes\OO_x\br=\ell\bl(\f\doubledual/\f
)\otimes\OO_x\br$\ for all $x\in X$. Then the image $H_i(V(\f))\util$\
is contained in $H_i(V_0)\util$.
\endpro

Now we prove the following proposition that is equivalent to
the surjectivity of (4.9).

\pro{Proposition 4.6}
Let $H_0$, $H_0\in\c\ssub\nsq^+$\ be fixed. Then there is an $N$\
such that for any $d\geq N$\ and $(H_0,I,d)$-suitable $H\in\c$,
the homomorphism
$$r(d)_i: H_i(\sox\mm)\lra H_i(\mm)
$$
is an isomorphism for $i\leq 2$.
\endpro

\proof
The statement for $i=0$\ follows from [GL,OG1]. The proof that
$r(d)_1$\ is an isomorphism is similar and easier than that of $r(d)_2$,
which we will prove now.

Let $N$\ be chosen so that all requirement of
$N$\ in the previous results have been met and the proposition
is true for $r(d)_1$.
The surjectivity of $r(d)_2$\ will be proved by first establishing
$H_2(\s\mm)\util\sub H_2(\so\mm)\util$\ and then
$H_2(\so\mm)\util\sub H_2(\sox\mm)\util$. The difficult in showing
the first inclusion lies in the fact that $\s\mm$\ is not locally
irreducible along $\sto\mm$. Thus $H_2(\so\mm)\to H_2(\s\mm)$\
is not surjective. However, there images in $H_2(\mm)$\ actually
coincide.

Before we go any further, we need to introduce some subsets
of $\s\mm$\ that will help us understand the geometry of the
compliment of $\so\mm\sub\s\mm$.
The first space is
$\pfd$\ that is an algebraic space whose closed points are
pairs
$$\{\e\sub\f\}: \f\in\mmm,\, \e\in\s\mm\ \text{and}\ \f/\e\cong\CC_p\ \text{
for some}\ p\in X.
$$
$\pfd$\ admits projections $\pim$\ and $\pix$\
onto $\s\mm$\ and $X$\ respectively by sending $\{\e\sub\f\}$\
to $\{\e\}\in\s\mm$\ and $\supp(\e/\f)$\ respectively. Note that
$\pfd$\ is irreducible and $\pim$\
is one-to-one restricting to $(\pim)\upmo(\soo\mm)$. The second space
is $\ptfd$\ consisting of tuples
$$\{\e\sub\fo\sub\ft\}: \f_2\in\mmmm, \fo/\e\cong \CC_{p_1}\
\text{and}\ \ft/\fo\cong\CC_{p_2}\ \text{for some}\ p_1, p_2\in X.
$$
We let $\pito\mh\ptfd\to\pfd$\
be the map sending $\{\e\sub\fo\sub\ft\}$\ to $\{\e\sub\fo\}$.
Let $\pitm(=\pim\circ\pito)$,
$\pixo$\ and $\pixt$\ be the projections from $\ptfd$\ to $\s\mm$\ and $X$\
respectively by sending $\{\e\sub\fo
\sub\ft\}$\ to $\e$, $p_1$\ and $p_2$\ respectively.

Now let $W_1\sub\st\mm-\sto\mm$\ and let $W_2\sub\s_3\mm$\ be
subsets consisting of points $w$'s so that $(\pitm)^{-1}(w)$\
are single point sets. we claim that $\dim \st\mm-(W_1\cup\sto\mm)
\leq \dim\mm-4$\ and $W_2\ne\emptyset$.
Indeed, let $\e\in\s\mm$\ be such that $\len(\e\doubledual/\e)=k\geq 2$.
Then $(\pitm)\upmo(\e)$\ is a point if and only if there
is a unique filtration $\t_0\sub\t_2\sub \t_2=\e\doubledual/\e$\
such that $\len(\t_j/\t_{j-1})=1$, $j=1,2$.
When $\len(\e\doubledual/\e)=2$\
and $\e\not\in\sto\mm$, then the uniqueness of the above filtration
is equivalent to $\bl\e\doubledual/\e\br\otimes\CC_p=\CC_p$\ for some
$p\in X$. From this description, it is easy to see that
$W_1$\ is dense in $\st\mm-\sto\mm$. Therefore,
$\dim \st\mm-(\sto\mm\cup W_1)\leq \dim\mm-4$.
To show $W_2\ne\emptyset$, one notices that any sheaf
that has the form $\OO_p\oplus (z_1,z_2^3)\OO_p$\ at $p\in X$, where
$(z_1,z_2)$\ is a local coordinate of $p\in X$, and locally free elsewhere
belongs to $W_2$. Such sheaves do exist in $\s_3\mm$. Finally, we
note that $\s_3\mm$\ is irreducible because its generic
points are kernel of $\f\to \oplus_{i=1}^3\CC_{p_i}$\
with $\f$\ locally free and $p_1,p_2,p_3\in X$\ distinct.

Now let $\xi\in H_2(\s\mm)$\ be any element. Since $\mm$\ is smooth and
$\s\mm$\ is Cartier, by lemma 1.11 and 1.12, we can find a representative
$f\mh \sig\to\s\mm$\ of a multiple of $\xi$\ such that $\sig$\
is a Riemann surface and

\noindent
(4.10):\ $f(\sig)\sub\so\mm\cup\sto\mm\cup W_1\cup W_2$;

\noindent
(4.11):\ $f\upmo(\sto\mm\cup W_1\cup W_2)$\ is an
(at most real 1-dimensional) stratified set.

Now let $C$\ be the closure of $f\upmo(\sto\mm)\sub\sig$\ and
let $\eta\mh S\to \sig$\ be the Riemann surface with boundary
obtained by cutting $\sig$\ along $C$. Because $\pim\mh\pfd\to\s\mm$\
is one-to-one over $f(\sig)-\sto\mm$, $f$\ lifts to a
unique $g\pri\mh S-\d S\to\pfd$\ and further, because $\pim$\ is finite
over $f(\sig)$, $g\pri$\ extends to $g\mh S\to \pfd$.
By the choice of $W_1$\ and $W_2$,
the preimage of any $z\in g(\d S)$\ of the map $\pito\mh
\ptfd\to\pfd$\ is a one point set. Therefore, if we let $\Gamma=\d S$,
the restriction of $g$\ to $\Gamma$\ lifts to
$\tilg\mh\Gamma\to\ptfd$, which gives rise to the following
commutative diagram of maps:
$$\CD
\Gamma @> \sub >> S @>\eta>> \sig\\
@V{\tilde g}VV @V g VV @V f VV\\
\ptfd @> \pito >> \pfd @> \pi^1_M >> \s\mm.
\endCD
$$

\pro{Lemma 4.7}
$\tilg\mh \Gamma\to\ptfd$\ represents the trivial class in $H_1(\ptfd)$.
\endpro

\proof
As chains,
$$\d\bl\pix\circ g(S)\br=\pix\circ g(\d S)=\pix\circ g(\Gamma).
$$
Thus $\pix\circ g(\Gamma)$\ represents the trivial class in $H_1(X)$.
Because $\pixo\circ\tilg =\pix\circ g$\ as maps from $\Gamma$\ to $X$.
$[\pixo\circ \tilde g(\Gamma)]=0$\
in $H_1(X)$\ as well.
We claim $\pixt\circ\tilg(\Gamma)$\ also represents the
trivial class in $H_1(X)$. Indeed, let $I_1,\cdots,I_k$\
(each is homeomorphic to [0,1]) be
segments of $C\sub S$\ with fixed orientations.
Then $\d S$\ can be divided into segments $I_1^{\pm},\cdots,I_k^{\pm}$\
such that $\eta(I_i^{\pm})=\pm I_i$\ as (oriented) chains and
$\d S=\sum_{i=1}^k(I_i^++I_i^-)$. Next, we let $\Lam_1$\
be those $I_i$\ such that $\text{Im}\,g(I_i^+)=\text{Im}\,g(I_i^-)$\
(which means that we do not need to cut $\Sig$\ along $I_i$\ in order to
get a lift $g$) and
let $\Lam_2$\ be the remainder $I_i$'s. Note that for $i\in \Lam_1$,
$g(I_i^+)=-g(I_i^-)$.
Then we have the following identities of chains:
$$\pixo\circ \tilde g(I_i^+)+\pixo\circ \tilde g(I_i^-)=
\pixt\circ \tilde g(I_i^+)+\pixt\circ \tilde g(I_i^-)=0,\ I_i\in\Lam_1;
$$
$$\pixo\circ g(I_i^{\pm})=-\pixt\circ \tilde g(I_i^{\mp}),\ I_i\in\Lam_2.
$$
Therefore, we have
$$[\pixt\circ \tilde g(\Gamma)]=-[\pixo\circ \tilde g(\Gamma)]=0
\in H_1(X).
\tag 4.12
$$

Now we consider the class $[\tilg(\Gamma)]\in H_1(\ptfd)$.
Let $p\mh\ptfd\to\mmmm$\ be the map sending $\{\e\sub\fo\sub\ft\}$\
to $\ft$. Because $\mmmm$\ is smooth and $\ptfd$\ is locally
irreducible, $\tilg(\Gamma)$\ is homotopic to a $\Gamma\pri\sub
\ptfd$\ such that $p(\Gamma\pri)\sub\mmmm^0$.
Let $\x=p\upmo(\mmmm^0)$. It is easy to check that fibers of $\x$\ over
$\mmmm^0\times X\times X$\ (via $p\times\pixo\times\pixt$) are connected and
have
trivial first homology groups. By Leray spectral sequence, $H_1(\x)$\
is isomorphic to $H_1(\mmmm^0
\times X\times X)$. Because of (4.12), for any pair of distinct points
$(x_1,x_2)\in X\times X$, there is a 1-cycle
$\Gamma''$\ contained in $\st^{x_1x_2}\mm=(\pixo\times\pixt)\upmo(x_1,x_2)$\
such that $\Gamma''$\ is homologous to $\tilde g(\Gamma)$\ in
$\ptfd$. Let the homology be given by $\sig\pri$. Namely,
$\Sig\pri\sub\ptfd$\ such that
$\d\sig\pri=\Gamma''-\tilde g(\Gamma)$. By lemma 1.11 and
1.12, we can choose $\sig\pri$\ such that
$\sig\pri\sub(\pitm)\upmo(\st\mm\cup W_2)$.
Finally, because $\pitm(\Gamma'')\sub \s_2^{x_1,x_2}\mm$\ is
the boundary of $\pitm(\sig\pri)\sub\s\mm$,
$[\pitm(\Gamma'')]\in H_1(\st^{x_1x_2}\mm)$\ is contained
in the kernel of $H_1(\st^{x_1x_2}\mm)\to H_1(\mm)$.
By the induction hypothesis,
$[\pitm(\Gamma'')]=0\in H_1(\st^{x_1x_2}\mm)$\ and then
$[\tilg(\Gamma)]=0\in H_1(\ptfd)$. This completes the proof  of the lemma.
\qed

\noindent
{\sl Continuation of the proof of proposition 4.6}.
{}From the proof, we see that we can find a 2-chain $\sig\pri
\sub\ptfd$\ so that $\d\sig\pri=\tilg(\Gamma)$\ and
$\pitm(\sig\pri)-\st\mm$\ is a discrete point set contained in $W_2$.
In summary, for any $\xi\in H_2(\s\mm)$,
we can find a closed stratifiable set $S\sub\pfd$\ such that
$\pim(S)$\ is a representative of a multiple of $\xi\in H_2(\mm)$\
and the boundary $\Gamma=\d S\sub \pfd$\ can be lifted to $\ptfd$, say
$\tilde \Gamma\to\ptfd$, with $[\tilde \Gamma]=0
\in H_1(\ptfd)$. Hence we can find a 2-chain $S\pri\sub\ptfd$\ such that
$$\d S\pri=-\tilde \Gamma,\,
\pitm(S\pri)\sub \st\mm\cup W_2\ \text{and}\
\pitm(S\pri)-\st\mm\sub W_2\ \text{is discrete}.
\tag 4.13
$$
Let $T_1=S\cup \pi^2_1(S\pri)\sub\pfd$\ and $T_2=\pitm(S\pri)
\sub\s\mm-\so\mm$. Then $u_1=[T_1]
\in H_1(\pfd)$\ and $u_2=[T_2]\in H_2(\s\mm-\so\mm)$\
has the property that their images $\bar u_1$\ and $\bar u_2$\
in $H_2(\mm)$\ provide a decomposition
$\xi=\bar u_1+\bar u_2$. Finally, by the construction,
$T_1$\ is contained in $(\pim)\upmo\bl\so\mm\cup\st\mm\cup W_2\br$\
and $T_2$\ is contained in $\st\mm\cup W_2$.

Our next step is to prove that $\bar u_1\in H_2(\soo\mm)\util$\
and $\bar u_2\in H_2(\sto\mm)\util$.
We first prove $\bar u_1\in H_2(\soo\mm)\util$.
We take the representative $T_1\sub\pfd$. By perturb $T_1$\
if necessary, we can assume $T_1-(\pim)\upmo(\so\mm)$\
is discrete whose image in $\s\mm$\ lies in $\st\mm\cup W_2$.
Let $\{p_1,\cdots,p_k\}\sub T_1$\ be these points and let
$\f_i$\ be the sheaf
in $\mmm$\ that corresponds to the image of $p_i$\ under the projection
$\pom\mh\pfd\to\mmm$\ by sending $\{\e\sub\f\}$\ to $\f$.
We let $V(\f)\sub\mmm$\ be the closed set defined in lemma 4.5.
Then there is an analytic deformation retract neighborhood
$U_i$\ of $V(\f)$\ such that $H_1(U_i-\s\mmm)=0$.
Because fibers of $\pfd\to\mmm$\ are $\PP^1$\ over $\mmmo$,
the homomorphism
$$H_2\bl(p_M)\upmo(\mmm^0)\br\oplus \oplu{k}{i=1} H_2(V(\f_i))\lra
H_2\Bigl((p_M)\upmo\bl\mmm^0\bigcup\cup_{i=1}^k V(\f_i)\br\Bigr)
$$
is surjective by Mayer-Vietoris sequence. However, we know that
$H_2(V(\f_i))\util\sub H_2(V_0)\util\sub H_2(\so\mm)\util$\ (lemma 4.5).
Therefore, $\bar u_1\in H_2(\so\mm)\util$.

Next, we prove $\bar u_2\in H_2(\st^0\mm)\util$. Because $\bar u_2$\
is the image of $[\pitm(S\pri)]$\ with $S\pri$\ satisfies (4.13) and because
$\st\mm$\ is locally irreducible, similar to
the previous reasoning, $\bar u_2$\ is contained in the linear span of
$H_2\bl\sto\mm\br\util$\ and $H_2\bl V(\f)\br\util$\ for all
possible $\f$\ in $W_2$\ and $\st\mm-\sto\mm$.
But by lemma 4.5, all these $H_2(V(\f))\util$\ are contained in $H_2(
\sto\mm)\util$. Hence, $\bar u_2\in H_2(\sto\mm)\util$.

It remains to show that
$$H_2(\so\mm)\util\sub H_2(\sox\mm)\util \quad \text{and}\quad
H_2(\sto\mm)\util\sub H_2(\so^x\mm)\util.
\tag 4.14
$$
We will prove the first inclusion and leave the proof of the second inclusion
to the readers. Since $\so\mm$\ is a $\PP^1$-bundle
over $X\times \mmm^0$\ with projection $\pi$, by Kunneth decomposition
$$\align
H_2\bl\so\mm\br & =\pi_{\ast}\upmo\bl H_2(X\times\mmmo)\br\\
&=\pi_{\ast}\upmo\bl H_2(\mmmo)\br
+\pi_{\ast}\upmo\bl H_1(\mmmo)\otimes H_1(X)\br+
\pi_{\ast}\upmo\bl H_2(X)\br,\\
\endalign
$$
where $\pi_{\ast}\mh H_2(\so\mm)\to H_2(X\times\mmmo)$.
Let $\pi_{\ast}\upmo\bl\cdot\br\util$\ be the image in $H_2(\mm)$\
of the respective space. We claim that all of them are contained in
$H_2(\sox\mm)\util$. First, $\pi_{\ast}\upmo\bl H_2(\mmmo)\util\br
=H_2\bl\sox\mm\br\util$\ by definition. Next, we fix a ball $B\sub X$\
containing $x$\ and a compact $S\sub\mmmo$\ such that
$$\pi_{\ast}\upmo\bl H_1(\mmmo)\otimes H_1(X)\br
\sub H_2(\w_{S}),
\tag 4.15
$$
where $\w_{S}$\ is the set of $\e\in\so\mm$\
such that for some $\f\in S$, $\e\sub\f$\ and
$\supp(\f/\e)\cap B=\emptyset$.
Because $H_1(\sox\mmm)\to H_1(\mmm)$\ is an isomorphism, by
induction hypothesis we can replace $S$\ by a set $S\pri
\sub\sox\mmm$\ and still have (4.15) with $S$\ replaced by $S\pri$.
On the other hand, $H_2(\w_{
S\pri})\util$\ is certainly contained in $H_2(\sox\mm)\util$.
Therefore, $\pi_{\ast}\upmo\bl H_1(\mmmo)\otimes H_1(X)\br\util
\sub H_2(\sox\mm)\util$. The inclusion $\pi_{\ast}\upmo\bl H_2(X)\br\util
\sub H_2(\sox\mm)\util$\ can be proved similarly. This completes
the proof of proposition 4.6.
\qed

Now we are ready to prove the main theorems. Let $H_0$\ be any ample
divisor, let $H_0\in\c\ssub\nsq^+$\ be a precompact neighborhood of
$H_0\in\nsq^+$\ and let $N$\ be the constant given by proposition 4.6 and
the preceding propositions. Then for any $H\in\c$\ and $d\geq N$,
$H_i(\mhid^0)$\ is isomorphic to $H_i(\MM_{H_0}(I,d)^0)$\ ($i\leq 2$\ here and
in the later discussion). However, for $d\geq N$\ and
$(H_0,I,d)$-suitable $H\in\c$, by theorem 4.1
$$H_i\bl\mhid^0,\QQ\br\lra H_i\bl\MM_H(I,d+1)^0,\QQ\br
$$
is surjective. Therefore,
$$H_i\bl\MM_{H_0}(I,d)^0,\QQ\br \lra H_i\bl\MM_{H_0}(I,d+1)^0,\QQ\br
\tag 4.16
$$
is surjective for all $d\geq N$. Since $H_i(\MM_{H_0}(I,d)^0)$\ are
linear spaces, the system (4.16) has to stabilize at some stage.
Namely, for some $N_1\geq N$, (4.16) is an isomorphism for all $d\geq N_1$.
Further, combined with the work of [Ta] (see (0.3) and (0.4)), we can
find an $N_2$\ so that for $d\geq N_2$,
$$H_i\bl\MM_{H_0}(I,d)^0,\QQ\br\cong H_i\bl\Cal B(P_d)\sta,\QQ\br.
$$
This proves theorem 0.1. For theorem 0.2, we simply apply
the above isomorphism to the fact
$$\dim H_1\bl\Cal B(P_d)\sta\br=b_1\ \text{and}\
\dim H_2\bl\Cal B(P_d)\sta\br=b_2+{1\over 2}b_1(b_1-1),
$$
where $b_i=\dim H_i(X)$\ (see page 181-182 of [DK]\footnote{Since SU(2)
or SO(3) is a rational three sphere, the first two rational homology groups
of $\tilde\Cal B\sta$\ and $\Cal B\sta$\ coincide.}).
Finally, by the proof of lemma 4.3, $h_1(\MM_{H_0}(I,d))=
h_1(\MM_{H_0}(I,d)^0)$, and $h_2 (\MM_{H_0}(I,d))=
h_2 (\MM_{H_0}(I,d)^0)+1$\ because $\s\MM_{H_0}(I,d)\sub
\MM_{H_0}(I,d)$\ is an irreducible Cartier divisor. This proves theorem 0.3.

In the remainder of this section, we shall prove lemma 4.5 promised
earlier. First, we will introduce set similar to $\Cal Z_2$\
that is a desingularization of $V(\f)$\ as topological space. Let
$\e_0$\ be any sheaf and let $x\in X$\ be fixed. We define $R_i$\ be the set
of all filtrations
$$\{\e_i\sub\e_{i-1}\sub\cdots\sub\e_0\}
\tag 4.17
$$
such that $\e_j/\e_{j-1}\cong\CC_x$. Let $\pi_{ij}\mh R_i\to R_j$, $i\geq j$,
be the map sending (4.17) to the subfiltration
$\{\e_j\sub\cdots\sub\e_0\}$\ and let
$p_i$\ be the map sending (4.17) to $\e_i$. Obviously, $R_1\cong \PP^1$\ and
for any $z\in R_1$, $\pi_{21}\upmo(z)\cong\PP^2$. Thus $R_2$\ is a
$\PP^2$-bundle
over $\PP^1$. Similarly, because for each $\{\e_2\sub\e_1\sub\e_0\}\in R_2$\
the tensor product
$\e_2\otimes\CC_x=\CC_x^{\oplus 4}$, the fibers of $\pi_{32}\mh R_3\to R_2$\
are isomorphic to $\PP^3$. Hence $R_3$\ is a $\PP^3$-bundle over $R_2$\
as topological spaces. Therefore, $R_2$\ and $R_3$\ are simply connected,
$H_2(R_2)=\QQ^{\oplus 2}$\ and $H_2(R_3)=\QQ^{\oplus 3}$.

Now we prove lemma 4.5. We first study the case where $\f\in\st\mm$\
and the support of $\f\doubledual/\f$\ is $x$. Let $\e_0=\f\doubledual$\
and let $R_2$\ be the set introduced based on $\e_0$\ and $x$. Then $p_2$\
maps $R_2$\ onto $V(\f)$. Let $w\in V(\f)$\ be the sheaf that is the kernel of
$\e_0\to\CC^{\oplus 2}_x$. Then $p_2$\ is one-to-one away from
$Z=p_2\upmo(w)$\ and $Z\cong P^1$. Hence $H_1(R_2)\to H_1(V(\f))$\ is
an isomorphism and $p_{2\ast}\mh H_2(R_2)\to H_2(V(\f))$\ is surjective
whose kernel contains $[Z]$.
Thus $\dim H_2(V(\f))\leq 1$. On the other hand,
$H_2(V_0)\util$\ is one dimensional and is obviously contained in
$H_2(V(\f))\util$.
Therefore, $H_2(V(\f))\util=H_2(V_0)\util$.
The case $\f\in\sto\mm$\ can be checked similarly.
This completes the proof of lemma 4.5 for $j=2$.

Now we study the case where $\f\in\s_3\mm$. We will deal with the case
where $\supp(\f\doubledual/\f)=\{x\}$\ and leave the other cases to the
readers.
Let $\e_0=\f\doubledual$\ and let $R_3$\ be as before. Then $p_3\mh R_3\to
V(\f)$\
is surjective that induces surjective homomorphisms
$H_1(R_3)\to H_1(V(\f))$\ and
$p_{3\ast}\mh H_2(R_3)\to H_2(V(\f))$. Since $H_2(R_3)=\QQ^{\oplus 3}$, as
before to prove the lemma it suffices to show that $\ker(p_{3\ast})$\ is
two dimensional. We will prove this by construct two
generators of $\ker(p_{2\ast})$\ explicitly. First, let $\e_2\pri$\
be the kernel of $\e_0\to\CC_x^{\oplus 2}$\ and let
$\e_3\pri\sub\e_2\pri$\ be any sheaf with $\e_2\pri/\e_3\pri\cong\CC_x$.
We define $Z_1$\ be the subset of $R_3$\ consists of filtrations of type
$$\{\e_3\pri\sub\e_2\pri\sub\e_1\sub\e_0\}.
$$
$Z_1$\ is isomorphic to $\PP^1$\ and $[Z_1]\ne0\in H_2(R_3)$\ belongs to
$\ker(p_{3\ast})$. Secondly, we let $w\in V(\f)$\ be a sheaf that is isomorphic
to
$\OO_x\oplus (u_1^2,u_1u_2,u_2^2)\OO_x$\ at $x$,
where $(u_1,u_2)$\ is an analytic coordinate of $x\in X$. Then $Z_2=p_3\upmo(
w)$\ is again a projective line that generates a nontrivial kernel of
$p_{3\ast}$. Since $\pi_{31}(Z_2)$\ is a point and
$\pi_{31}(Z_1)=R_1$, $[Z_1]$\ and $[Z_2]$\
generates a two dimensional subspace in $H_2(R_3)$. Therefore
$\dim H_2(V(\f))\leq 1$. Finally, because
$\QQ=H_2(V_0)\util\sub H_2(V(\f))\util$, they must be identical.
This proves the lemma 4.5.

\end